\documentclass{amsart}

\newcommand{\cal}{\mathcal}

\newtheorem{teo}{Theorem}
\newtheorem{lem}{Lemma}
\newtheorem{pro}{Proposition}
\newtheorem{cor}{Corollary}

\newtheorem{defi}{Definition}

\theoremstyle{definition}

\newtheorem{examp}{Example}

\def\sp{\hspace{0.2cm}}
\def\qed{\begin{flushright} $\square$ \end{flushright}}
\def\qee{\begin{flushright} $\Diamond$ \end{flushright}}
\def\ov{\overline}

\begin{document}

\title[Quantum Stochastic Processes, Quantum IFS and Entropy]
{Quantum Stochastic Processes, Quantum Iterated Function Systems and Entropy}

%\author{A. Baraviera, C. F. Lardizabal,  A. O. Lopes, and M. Terra Cunha}

\author{A. Baraviera}
\address{}
\curraddr{}
\email{}
\thanks{}

\author{C. F. Lardizabal}
\address{}
\curraddr{}
\email{}
\thanks{Supported in part by CAPES and CNPq}

\author{A. O. Lopes}
\address{}
\curraddr{}
\email{}
\thanks{}

\author{M. Terra Cunha}
\address{}
\curraddr{}
\email{}
\thanks{}

\begin{abstract}
We describe some basic results for Quantum Stochastic Processes
and present some new results about a certain class of processes which are
associated to Quantum Iterated Function Systems (QIFS). We
discuss questions related to the Markov property and we present a definition of entropy which is induced by a QIFS. This definition is a natural generalization of
the Shannon-Kolmogorov entropy from Ergodic Theory.

\end{abstract}

\maketitle

\bigskip

{Paper to appear in S\~ao Paulo Journal of Mathematical Sciences
(2010)}

\bigskip

\section{Introduction}

We review and discuss some main properties of  Quantum Stochastic
Processes (see \cite{GZ} \cite{slom} \cite{Sr}) and present some new results about
a certain class of processes which are associated to a Quantum Iterated Function System (QIFS). The concept of QIFS  was introduced in the work \cite{lozinski}, and it is a natural object in Quantum Information Theory.

We also present a definition of entropy which
is suitable for the QIFS. This definition is a natural
generalization of the Shannon-Kolmogorov entropy of Ergodic
Theory. We describe a parallel between the classical Kolmogorov
entropy and the one we present here, which is different from the one seen in \cite{BLLT1}.

The present  definition of entropy is obtained by adapting the reasoning
described in \cite{GL}, \cite{lopes} and \cite{lopes_elismar} to
the setting we present in this work. The main idea is to define this concept via the
Ruelle operator and to avoid the use of partitions. Using this definition one can consider maximal pressure density states. This formulation can be seen as a mini-max problem (see \cite{GL} \cite{lopes} \cite{lopes_elismar}). In \cite{BLLT1} it is
described some applications of the pressure problem.

Section \ref{s_not} introduces basic notations and section \ref{qsp_slom} describes QSPs following \cite{slom}; section \ref{s_qifs} and \ref{ex_qifs} describes Quantum Iterated Function Systems, following \cite{lozinski}. Section \ref{pqeck} is a brief digression on the Chapman-Kolmogorov equation and probability amplitudes. Section \ref{s_piqifs} defines probabilities measures and quantum stochastic processes induced by QIFS. Section \ref{entr_artur} gives a definition of entropy induced by a QIFS and we make a few remarks on the variational problem of pressure.

Our work is  inspired by results presented in \cite {lozinski}
and \cite{wsbook}. We would like to thank these authors for
supplying us with the corresponding references. Some other references related to the topic described here are \cite{BLLT1} \cite{BLLT2} \cite{benatti} \cite{bcs} \cite{wsbook}.

This work is part of the thesis dissertation of C. F. Lardizabal in Prog. Pos-Grad. Mat. UFRGS (Brazil) \cite{lar}.

\section{Notations}\label{s_not}

We recall some basic notation which is used in Quantum Computing.
For a comprehensive introduction to the subject, see \cite{nich}.
Let $\mathcal{H}_N$ be a Hilbert space of finite dimension $N$. If
a quantum system is in a certain known state $\vert
\psi\rangle\in\mathcal{H}_N$, we say that the system is in a pure
state. Otherwise the system is in a mixed state. Each system
contains certain pure states, which are fixed when we define our
problem. Also, such states are normalized, so we have
$\langle\psi\vert\psi\rangle=1$. For any phase $\alpha$, we
identify the elements
$\vert\psi'\rangle=e^{i\alpha}\vert\psi\rangle$ and
$\vert\psi\rangle$, so we get the space of pure states, denoted by
$\mathcal{P}_N$. Topologically, it is the complex projective space
$\mathbb{CP}^{N-1}$ with the Fubini-Study metric, given by
$D_{FS}(\vert\phi\rangle,\vert\psi\rangle):=
\text{arccos}\vert\langle\phi\vert\psi\rangle\vert.$

A {\bf qubit} is a unit vector in a complex vector space of dimension 2
$$\vert\psi\rangle=\alpha\vert 0\rangle +\beta\vert 1\rangle,$$
where $\vert\alpha\vert^2+\vert\beta\vert^2=1$. We can rewrite such equation as
$$\vert\psi\rangle=e^{i\gamma}(cos\frac{\theta}{2}\vert 0\rangle + e^{i\phi}sin\frac{\theta}{2}\vert 1\rangle),$$
where $\theta, \phi, \gamma$ are real numbers. As we are in projective space, the factor $e^{i\gamma}$ can be ignored, so we can write
$$\vert\psi\rangle=cos\frac{\theta}{2}\vert 0\rangle + e^{i\phi}sin\frac{\theta}{2}\vert 1\rangle$$
The numbers $\theta$ and $\phi$ define a point on the unit sphere, the Bloch sphere, which gives us an easy way to visualize the state of a qubit.

\bigskip

Denote by $\rho^*$ the adjoint of $\rho:\mathcal{H}_N\to\mathcal{H}_N$. We say that $\rho:\mathcal{H}_N\to\mathcal{H}_N$ is hermitian if $\rho=\rho^*$. We say that a hermitian operator $P:\mathcal{H}_N\to \mathcal{H}_N$ is positive, denoting such fact by  $P\geq 0$, if
$\langle Pv,v\rangle\geq 0, \sp\forall v\in\mathcal{H}_N$.

\begin{defi}
A {\bf density operator} (or density matrix) is an operator $\rho$ acting on $\mathcal{H}_N$, with $\rho=\rho^*$, $\rho\geq 0$ and $tr\rho=1$. Denote by $\mathcal{M}_N$ the space of density operators.
\end{defi}
If $\vert\psi\rangle\in\mathcal{H}_N$ denote the associated projection by $\vert \psi\rangle\langle\psi\vert$. We denote by
$$\{\vert 0\rangle,\dots, \vert N-1\rangle\},$$
the canonical orthonormal basis for $\mathcal{H}_N$. A density operator $\rho$ can always be written as
\begin{equation}
\rho=\sum_{i=1}^k p_i\,\,\vert \psi_i\rangle\langle \psi_i\vert
\end{equation}
where the $p_i$ are positive numbers with $\sum_i p_i=1$, and $\psi_i$, $i=1,2,...N-1$, is an orthonormal basis.

A pure state is such that its associated density operator satisfies $tr(\rho^2)=1$; if a state is mixed, we have $tr(\rho^2)<1$. Also an operator is a density operator if and only if its trace equals 1 and if it is positive.

\section{A description of quantum stochastic process}\label{qsp_slom}

In this section the definitions and examples were taken from
\cite{slom}, where it is presented  a definition of quantum
stochastic process. We briefly describe some of the results
obtained in that work.

\begin{defi}
A state space is a pair $(V,K)$, where
\begin{enumerate}
\item $V$ is a real Banach space with norm $\Vert\cdot\Vert$.
\item $K$ is a closed cone in $V$.
\item If $u,v\in K$ then $\Vert u\Vert+\Vert v\Vert=\Vert u+v\Vert$
\item If $u\in V$ e $\epsilon >0$ then there are $u_1,u_2\in K$ such that $u=u_1-u_2$ and $\Vert u_1\Vert+\Vert u_2\Vert <\Vert u\Vert +\epsilon$.
\end{enumerate}
\end{defi}

\begin{defi} If $(V,K)$ is a state space then there is a unique positive linear functional  $\tau: V\to \mathbb{R}$ such that $\tau(u)=\Vert u\Vert$ if $u\in K$, and $\tau(u)\leq \Vert u\Vert$ if $u\in V$. We say that $u\in K$ is a {\bf state} if $\tau(u)=1$.
\end{defi}

\begin{examp} Let $\mathcal{H}$ be a finite dimensional Hilbert space and let $V$ be the space of hermitian operators in $\mathcal{H}$. Let $K$ be the set of positive operators in $V$. In this case we have $\tau(B)=tr(B)$ for all $B$ operator in $V$.

\end{examp}

\qee

\begin{defi}
A {\bf phase space} is a measurable space $(\Omega,\Sigma)$ where $\Omega$ represents the set of all possible results for a measurement and $\Sigma$ is a $\sigma$-algebra of subsets of $\Omega$.
\end{defi}

 Let $V^*$ be the dual space of $V$. We introduce a partial order on $V^*$ by defining $\phi\geq\psi$ if $\phi(u)\geq\psi(u)$, for all $u\in K$.

\begin{defi}
   An {\bf effect} is a mapping $\phi\in V^*$ such that $0\leq\phi\leq\tau$. We denote the space of effects by $\mathcal{E}\subset V^*$.
\end{defi}

\begin{defi}
  We say that $x:\Sigma\to\mathcal{E}$ is an {\bf observable} if $x$ is a measure taking values on the space of effects, such that $x(\Omega)=\tau$.
\end{defi}

If $E\in\Sigma$, $u\in K$ and $\tau(u)=1$ then $x(E)u$ can be interpreted as the probability that the result of the measurement of the physical quantity represented by $x$, prepared in the state $u$, belongs to the set $E$. In the case of quantum mechanics in Hilbert space, effects can be identified with bounded operators $A$ such that $0\leq A\leq 1$ by the formula $\phi_A(W)=tr(AW)$.

\begin{defi}
 An {\bf operation} is a positive linear operator $T:V\to V$ satisfying $0\leq\tau(Tu)\leq \tau(u)$ for all $u\in K$. The space of operations will be denoted by $\mathcal{O}$.
 \end{defi}

\begin{defi}
An {\bf operator valued measure}, or an OVM on a phase space is a map $\mathcal{I}:\Sigma\to\mathcal{O}$ such that if $\{E_n\}$ is a sequence of disjoint sets in $\Sigma$, then $\mathcal{I}(\cup E_n)=\sum\mathcal{I}(E_n)$.
\end{defi}

\bigskip
\begin{defi}
 Let $\mathcal{I}:\Sigma\to\mathcal{O}$ be an OVM, then we say that $\mathcal{I}$ is an
 {\bf instrument} if
\begin{equation}\label{instrumento_def}
\tau(\mathcal{I}(\Omega)u)=\tau(u), \forall u\in V.
\end{equation}
\end{defi}

We interpret such notion in the following way. Let $\mathcal{I}$ be an instrument, $E\in\Sigma$, $u\in K$. If $u$ is the state of the system before the measurement and if  $\mathcal{I}$ determines a value in $E$ then the resulting state is given by
\begin{equation}
\frac{\mathcal{I}(E)u}{\tau(\mathcal{I}(E)u)}
\end{equation}
Note that for each instrument $\mathcal{I}$, there is a unique observable $x_{\mathcal{I}}:\Sigma\to\mathcal{E}$ such that $\tau(\mathcal{I}(E)u)=x_{\mathcal{I}}(E)u$, $E\in\Sigma$, $u\in K$. Also, it is possible that two instruments correspond to the same observable \cite{slom}.

\bigskip

The following are examples of instruments:

\bigskip

\begin{examp}\label{exe1} Let $\mathcal{H}$ be a Hilbert space, and let $\mathcal{F}(\mathcal{H})$  be the space of hermitian operators $A$ in $\mathcal{H}$ such that
$$\sum_{k\in\mathbb{N}}\langle e_k, Ae_k\rangle<\infty$$
and have the same value in any orthonormal base $\{e_k\}_{k\in\mathbb{N}}$ for $\mathcal{H}$.
Let $\Omega=\{1,\dots, N\}$, or $\Omega=\mathbb{N}$, let $\{P_i\}_{i\in\Omega}$ be a family of orthogonal projections such that $\sum_i P_i=I$. Define
$$\mathcal{I}:\Sigma\to\mathcal{O}$$
$$x_{\mathcal{I}}:\Sigma\to\mathcal{E}$$
as
\begin{equation}
\mathcal{I}(E)\rho:=\sum_{i\in E} P_i\rho P_i,
\end{equation}
\begin{equation}
x_\mathcal{I}(E)\rho:=\sum_{i\in E}\tau(P_i\rho),
\end{equation}
for all $E\subset\Omega$ and $\rho\in \mathcal{F}(\mathcal{H})$.
\end{examp}

\qee

\begin{examp}\label{exe2} Let $\mathcal{H}$ be a Hilbert space, $\Omega$ a topological space, $\Sigma$  a $\sigma$-algebra for $\Omega$ and $m$ a measure on $(\Omega,\Sigma)$. Let $\{P_a\}_{a\in\Omega}$ be a family of projections on $\mathcal{H}$, such that the mapping
$a\to P_a$ is strongly continuous and $\int_\Omega P_a dm(a)=I$. Then define
$$\mathcal{I}:\Sigma\to\mathcal{O}$$
$$x_{\mathcal{I}}:\Sigma\to\mathcal{E}$$
as
\begin{equation}
\mathcal{I}(E)\rho:=\int_E P_a\rho P_a dm(a)
\end{equation}
\begin{equation}
x_{\mathcal{I}}(E)\rho:=\int_E \tau(P_a\rho) dm(a),
\end{equation}
for all $E\subset\Omega$ e $\rho\in \mathcal{F}(\mathcal{H})$.
\end{examp}

\qee

\begin{examp}\label{exe3} Let $X$ be a locally compact Hausdorff space, $V$ the space of the countably additive functions on the Borel $\sigma$-algebra $\mathcal{B}(X)$ for $X$ endowed with the norm of total variation. Let $K$ be the set of nonnegative measures on $V$. Let $(\Omega,\Sigma)=(X,\mathcal{B}(X))$. Then
\begin{equation}
\mathcal{I}(E)\mu(A)=\mu(A\cap E),
\end{equation}
for $\mu\in V$, $A, E\in \Sigma$ is an instrument, called the {\bf sharp classical measurement} and the corresponding observable is
\begin{equation}
x_{\mathcal{I}}(E)\mu=\mu(E)
\end{equation}
\end{examp}

\qee

\begin{defi}
    Following \cite{slom}, a {\bf Quantum Stochastic Process}, QSP, is an arbitrary family of instruments $\{\mathcal{I}_t\}_{t\in\mathcal{J}}$. Let $\mathcal{J}=\mathbb{Z}$ or $\mathcal{J}=\mathbb{R}$ for discrete or continuous time, respectively.
\end{defi}

\bigskip
 The {\bf finite dimensional distributions} of the process are measures $\mu^{u}_{t_0,\dots,t_{n-1}}$ defined in $(\Omega^n,B(\Omega^n))$ as being the natural extensions of the functions given by
\begin{equation}\label{dff_1}
\mu^{u}_{t_0,\dots,t_{n-1}}(E_0\times\cdots\times E_{n-1})=\tau((\mathcal{I}_{t_{n-1}}(E_{n-1})\circ\mathcal{I}_{t_{n-2}}(E_{n-2})\circ\cdots\circ\mathcal{I}_{t_{0}}(E_{0}))u) \end{equation}
where $n\in\mathbb{N}$, $t_0\leq\cdots\leq t_{n-1}$, $t_i\in\mathcal{J}$, $u\in V$ and $E_0,\dots ,E_{n-1}\in\Sigma$. The meaning of such expression is the following: $\mu^{u}_{t_0,\dots,t_{n-1}}(E_0\times\cdots\times E_{n-1})$ is the joint probability that successive measurements of the system by the instruments $\mathcal{I}_0, \dots ,\mathcal{I}_{n-1}$ in the moments $t_0,\dots ,t_{n-1}$ produce values in $E_0,\dots , E_{n-1}$, when the pre-measurement state is $u$.

\bigskip

A {\bf probability transition} is a function $P:\Omega\times\Sigma\to\mathbb{R}$ such that $P(\cdot,E)$ is measurable for all $E\in\Sigma$ and $P(x,\cdot)$ is a probability measure for all $x\in\Omega$.

\begin{defi}
We say that a QSP is {\bf Markov} if there exists a family of probability transitions
$\{P_{s,t}\}_{s<t}$ such that
$$\mu^{u}_{t_0,\dots,t_{n-1}}(E_0\times\cdots\times E_{n-1})$$
\begin{equation}
=\int_{E_0}\int_{E_1}\cdots \int_{E_n} P_{t_{n-1},t_n}(y_{n-1},dy_n)\cdots P_{t_0,t_1}(y_0,dy_1)\mu_{t_0}^u(dy_0)
\end{equation}
 for all $t_0<\dots <t_n, t_i\in\mathcal{J}$, $u\in V$, $E_0,\dots, E_n\in\Sigma$.  A Markov QSP is {\bf homogeneous} if the probability transitions $P_{s,t}$ depend only on the difference $t-s$.
\end{defi}

\bigskip
{\bf Remark} In contrast with the classic theory of stochastic processes, the probability transitions of a Markov QSP do not satisfy in general the Chapman-Kolmogorov equation.
\qee

\begin{defi}
Let $\mathcal{I}$ be an instrument. Assume that between the measurements the system evolves and its evolution is described by a group $\{T_t\}_{t\in\mathcal{J}}$ of isometric automorphisms of $V$. Then define the QSP $\{\mathcal{I}_t\}_{t\in\mathcal{J}}$, where
\begin{equation}
\mathcal{I}_t(E)=T_t^{-1}\circ \mathcal{I}(E)\circ T_t
\end{equation}
is called a {\bf transformed instrument}. For simplicity, we can choose $\mathcal{J}=\mathbb{Z}$ so $T_n=T^n$ and we denote such process by $\mathcal{C}(T,\mathcal{I})$.
\end{defi}

Now we show an example of a Markov QSP.

\begin{examp}
Let $\mathcal{I}$ be the instrument given in example \ref{exe3} and let $\Theta: X\to X$ be a measurable map. Then $\Theta$ generates an automorphism $T_{\Theta}: V\to V$ by
\begin{equation}
T_{\Theta}(\mu)(A)=\mu(\Theta^{-1}(A)), \sp \mu\in V, A\in B(X)
\end{equation}
Then we can show that $\mathcal{C}(T_{\Theta},\mathcal{I})$ is a homogeneous Markov QSP and its transition probability is given by
\begin{equation}
P(x,E)=\mathcal{X}_E(\Theta x),\sp x\in X,\sp E\in B(X)
\end{equation}

\end{examp}

\qee

\section{Quantum IFS}\label{s_qifs}

This section follows \cite{lozinski}. We begin with a few definitions.

\begin{defi}
Let $G_i:\mathcal{M}_N\to\mathcal{M}_N$, $p_i:\mathcal{M}_N\to[0,1]$, $i=1,\dots ,k$ and such that $\sum_i p_i(\rho)=1$. We call
\begin{equation}\label{qifs}
\mathcal{F}_N=\{\mathcal{M}_N, G_i, p_i:i=1,\dots, k\}
\end{equation}
a {\bf Quantum Iterated Function System} (QIFS).
\end{defi}

\begin{defi}
A QIFS is {\bf homogeneous} if $p_i$ and $G_ip_i$ are affine mappings, $i=1,\dots , k$.
\end{defi}

Suppose that the QIFS considered is such that there are $V_i$ and $W_i$ linear maps, $i=1,\dots, k$, with $\sum_{i=1}^k W_i^*W_i=I$ such that
\begin{equation}
G_i(\rho)=\frac{V_i\rho V_i^*}{tr(V_i\rho V_i^*)}
\end{equation}
and
\begin{equation}
p_i(\rho)=tr(W_i\rho W_i^*)
\end{equation}
Then we have that a QIFS is homogeneous if $V_i$=$W_i$, $i=1,\dots,k$. Now we can define a Markov operator
$P:\mathcal{M}(\mathcal{M}_N)\to\mathcal{M}(\mathcal{M}_N)$,
$$(P\mu)(B)=\sum_{i=1}^k\int_{G_i^{-1}(B)}p_i(\rho)d\mu(\rho),$$
where $\mathcal{M}(\mathcal{M}_N)$ denotes the space of probability measure over $\mathcal{M}_N$. We also define $\Lambda: \mathcal{M}_N\to\mathcal{M}_N$,
$$\Lambda(\rho):=\sum_{i=1}^k p_i(\rho)G_i(\rho)$$
If the QIFS considered is homogeneous then
\begin{equation}\label{lambda_um}
\Lambda(\rho)=\sum_i V_i\rho V_i^*
\end{equation}

We say that $\rho\in \mathcal{M}_N$ is the {\bf integral} of a mapping $f:\mathcal{M}_N\to \mathcal{M}_N$, denoted by
$$\rho:=\int_{\mathcal{M}_N} fd\mu$$
if
$$l(\rho)=\int_{\mathcal{M}_N}l\circ fd\mu,$$
for all $l\in \mathcal{M}_N^*$.

\begin{teo} A mixed state $\rho_0$ is $\Lambda$-invariant, if and only
if,
\begin{equation}\label{bari}
\rho_0=\int_{\mathcal{M}_N} \rho d\mu(\rho),
\end{equation}
for some $P$-invariant measure $\mu$.
\end{teo}
For the proof, see \cite{lozinski}, \cite{wsbook}.

\bigskip

In order to define hyperbolic QIFS, we have to specify a distance on the space of mixed states. The following are a few possibilities:
$$D_1(\rho_1,\rho_2)=\sqrt{tr[(\rho_1-\rho_2)^2]}$$
$$D_2(\rho_1,\rho_2)=tr\sqrt{(\rho_1-\rho_2)^2}$$
$$D_3(\rho_1,\rho_2)=\sqrt{2\{1-tr[(\rho_1^{1/2}\rho_2\rho_1^{1/2})^{1/2}]\}}$$
Such metrics generate the same topology on $\mathcal{M}_N$. Considering the space of mixed states with one of those metrics we can make the following definition.

\begin{defi}
  We say a QIFS is {\bf hyperbolic} if the quantum maps $G_i$ are contractions with respect to one of the distances on $\mathcal{M}_N$ and if the maps $p_i$ are H\"older-continuous and positive.
\end{defi}

\begin{pro}\cite{lozinski} \cite{wsbook} If a QIFS (\ref{qifs}) is homogeneous and hyperbolic then the associated Markov operator admits a unique invariant measure $\mu$. Such invariant measure determines a unique  $\Lambda$-invariant state $\rho\in\mathcal{M}_N$, given by (\ref{bari}).
\end{pro}

\section{Examples of QIFS}\label{ex_qifs}

\begin{examp}\label{exemp7} $\Omega=\mathcal{M}_N$, $k=2$, $p_1=p_2=1/2$, $G_1(\rho)=U_1\rho U_1^*$, $G_2(\rho)=U_2\rho U_2^*$. The normalized identity matrix $\rho_*=I/N$ is $\Lambda$-invariant, for any choice of unitary $U_1$ and $U_2$. Note that we can write
$$\rho_*=\int_{\mathcal{M}_N}\rho d\mu(\rho)$$
where the measure $\mu$, uniformly distributed over $\mathcal{P}_N$, is $P$-invariant.
\end{examp}

\qee

\begin{examp}\label{exemp8} Let $\Omega=\mathcal{M}_N$, $k=2$, $p_1=p_2=1/2$, $G_1(\rho)=(\rho+2\rho_1)/3$, $G_2(\rho)=(\rho+2\rho_2)/3$, where we choose the projectors  $\rho_1=\vert 1\rangle\langle 1\vert$ and $\rho_2=\vert 2\rangle\langle 2\vert$ so that they are orthogonal. Since $G_1$ and $G_2$ are contractions with Lipschitz constant equal to $1/3$, this QIFS is hyperbolic and so there is a unique invariant measure.
\end{examp}

\qee

Recall that a mapping $\Lambda$ is {\bf completely positive} (CP), if $\Lambda\otimes I$ is positive for any extension of the original Hilbert space $\mathcal{H}_N\to\mathcal{H}_N\otimes\mathcal{H}_E$. We know that every trace preserving CP map can be represented (in a nonunique way) in the Stinespring-Kraus form
$$\Lambda(\rho)=\sum_{j=1}^k V_j\rho V_j^*, \sp \sum_{j=1}^k V_j^* V_j=1,$$
where the $V_j$ are linear operators.  Besides, if $\sum_{j=1}^k V_j  V_j^*=I$ then $\Lambda(I/N)=I/N$ and $\Lambda$ will be called {\bf unital}. This is the case if each of the $V_j$ is normal, that is, if $V_j V_j^*=V_j^* V_j$. Note that by writing $G_i(\rho)=U_i\rho U_i^*$, we have that example \ref{exemp7} is contained in this class of  QIFS. We call such QIFS {\bf unitary}. For a unitary QIFS we have that $\rho_*$ is an invariant state for $\Lambda_U$ and also that $\delta_{\rho_*}$ is invariant for the Markov operator $P_U$ induced by this QIFS.

\bigskip
\begin{defi}
   We say that unitary matrices of same dimension are {\bf common block diagonal} if they are block diagonal in the same base and with the same blocks.
\end{defi}

The proof of the following lemma is presented in \cite{lozinski}.

\begin{pro}\label{blocolema}
Assume that $p_i$, $i=1,\dots ,k$ are strictly positive. The the maximally mixed state  $\rho_*$ is the unique invariant state for the operator $\Lambda_U$ if and only if the unitary operators $U_i$, $i=1,\dots, k$ are not common block diagonal.
\end{pro}

\begin{examp}\label{exemp9} Let $\Omega=\mathcal{P}_2$, $U_1=I$, $U_2=\sigma_1$, $U_3=\sigma_2$, $U_4=\sigma_3$, $p_1=1-p$, $p_2=p_3=p_4=p/3>0$, where
$\sigma_1, \sigma_2, \sigma_3$ are the Pauli matrices. Since such matrices are not common block diagonal the maximally mixed state $\rho_*$ is the unique invariant state for the mapping below, called a quantum depolarizing channel \cite{lozinski}:
$$\Lambda_U(\rho)=\sum p_iU_i\rho U_i^*=(1-p)\rho+\frac{p}{3}(\sigma_1\rho\sigma_1+\sigma_2\rho\sigma_2+\sigma_3\rho\sigma_3).$$
\end{examp}

\qee

\begin{examp}\label{exemp10} Let $\Omega=\mathcal{P}_2$, $p_1=1-p$, $p_2=p$,
$$U_1=exp(-iH_0T/\hbar),$$
$$U_2=exp(-\frac{i}{\hbar}(H_0T +\int_0^T V(t)dt))$$
where $V(t)=V(t+T)$. The maximally mixed state $\rho_*=I/2$ is an invariant state for the operator $\Lambda_U$ corresponding to this QIFS. For a generic perturbation $V$, matrices $U_1$ and $U_2$ are not common block diagonal so $\rho_*$ is the unique invariant state for $\Lambda_U$.
\end{examp}

\qee

\section{On certain probability and amplitude calculations}\label{pqeck}

We begin with a
brief digression on the Chapman-Kolmogorov equation. Let $X=\{
X_n\}$ be a sequence of measurable functions. Suppose that
$$P(X_{n+1}=j\vert X_n=i)=P(X_1=j\vert X_0=i)$$
for all $n,i,j$. Suppose that $X$ takes values on a finite set $S$. Define the matrix $P=(p_{ij})$ of order $\vert S\vert$, with entries
$$p_{ij}=P(X_{n+1}=j\vert X_n=i)$$
Define the matrix of $n$ transitions $P_n=(p_{ij}(n))$, where
$$p_{ij}(n)=P(X_{m+n}=j\vert X_m=i)$$
Also suppose that it is a Markov chain, that is
\begin{equation}
P(X_n=x_n\vert X_0=x_0,X_1=x_1,\dots, X_{n-1}=x_{n-1})=P(X_n=x_n\vert X_{n-1}=x_{n-1})
\end{equation}
for all $n\geq 1$, and $x_0,\dots, x_n\in S$.

By using the fact that for any events $A_1,A_2,A_3$, we have
\begin{equation}\label{p_bas_prob01}
P(A_1\cap A_2\vert A_3)=P(A_1\vert A_2\cap A_3)P(A_2\vert A_3)
\end{equation}
we can write
$$p_{ij}(m+n)=P(X_{m+n}=j\vert X_0=i)=\sum_k P(X_{m+n}=j,X_m=k\vert X_0=i)$$
\begin{equation}\label{pqpqpq1}
=\sum_k P(X_{m+n}=j\vert X_m=k)P(X_m=k\vert X_0=i)
\end{equation}
So
\begin{equation}\label{chap_kolm}
p_{ij}(m+n)=\sum_k p_{ik}(m)p_{kj}(n)
\end{equation}
which is the Chapman-Kolmogorov equation. We are interested in studying quantum stochastic processes and in obtaining an adequate definition to what we will call a Markov quantum stochastic process. First we recall that in the previous section we have presented a description \cite{slom} of Markov QSP in which the Chapman-Kolmogorov do not hold in general. This fact can be seen as the general rule for quantum processes (but see \cite{gudder_pap2} for different settings).

\bigskip

In algebraic terms, we can argue that the deduction of
(\ref{chap_kolm}) above is not valid for quantum processes because
of equation (\ref{p_bas_prob01}). Since we have to take in
consideration the interference between measurements, the problem
of understanding how probability measures work in a quantum
setting is a basic question. In quantum mechanics we could in
principle consider a probability space $(\Omega,\Lambda,\mu)$ such
as in classic measure theory. However, we have that $\Lambda$ is a
$\sigma$-algebra and $\mu$ is a measure on $\Lambda$  only when we
are restricted to a single measurement. When we perform several
measurements interference effects occur and so we are no longer
considering a problem on classic probability \cite{gudder_book}.
Results of more general nature are presented in \cite{gudder3}.

\bigskip

We can think that interference occurs because, in contrast to classic probability measures, which can
be quite arbitrary, quantum probability measures are obtained in a very specific way.
In quantum mechanics we have an amplitude function $a:\Omega\to\mathbb{C}$, and if $B\in\Lambda$, we define the {\bf amplitude} of $B$ as
\begin{equation}
A(B)=\sum_{\omega\in B} a(\omega)
\end{equation}
and we define the probability that $B$ occurs as
\begin{equation}
\mu(B)=\vert A(B)\vert^2
\end{equation}
Let us describe a few more details on this point. For more on the subject, see for instance \cite{gudder_book}. Let $\Omega$ be a nonempty set and let $a:\Omega\to\mathbb{C}$. We say that $\omega\in\Omega$ is a  {\bf sample point} and the map $a$ is a {\bf probability amplitude}, and $(\Omega,f)$  is called a {\bf quantum probability space}. A set $A\subset \Omega$ is {\bf summable} if $\sum_{\omega\in\Omega}\vert a(\omega)\vert^2<\infty$ and we denote the collection of summable sets by $\Sigma_0$. Now define $A:\Sigma_0\to\mathbb{C}$ as $A(\emptyset)=0$ and
\begin{equation}
A(B):=\sum_{\omega\in\Omega}a(\omega)
\end{equation}
We say that $A(B)$ is the {\bf amplitude} of $B$. Now define
\begin{equation}
A(B_1\vert B_2):=\frac{A(B_1\cap B_2)}{A(B_2)}
\end{equation}
if $A(B_2)\neq 0$ and equal to zero, otherwise. In the case that $A(B_2)\neq 0$, we have that  $A(\cdot\vert B_2)$ is a complex measure on $P(\Omega)$, with $A(\Omega\vert B_2)=1$. We say that $A(B_1\vert B_2)$  is the {\bf conditional amplitude} of $B_1$, given $B_2$. Note that $A(B)=0$ does not imply $A(B\cap C)=0$ \cite{gudder_book}. Because of that, formulas of the kind $A(B\cap C)=A(B)A(C\vert B)$ might not be true when $A(B)=0$. However, when the conditioning sets have a nonzero amplitude, we have the formula
\begin{equation}
A(B_1\cap \cdots \cap B_n)=A(B_1)A(B_2\vert B_1)A(B_3\vert B_1\cap B_2)\cdots A(B_n\vert B_1\cap\cdots\cap B_{n-1})
\end{equation}
which is the amplitude counterpart for equation (\ref{p_bas_prob01}). Define the matrix $A=(a_{ij})$, where
$a_{ij}=A(X_{n+1}=j\vert X_n=i)$.
Now suppose that the chain $\{X_n\}_{n\in \mathbb{N}}$ is {\bf quantum Markov}, that is,
\begin{equation}
A(X_n=x_n\vert X_0=x_0,X_1=x_1,\dots, X_{n-1}=x_{n-1})=A(X_n=x_n\vert X_{n-1}=x_{n-1})
\end{equation}
 for all $n\geq 1$, $x_0,\dots, x_n\in S$. So in a way which is similar to what we did for probabilities, define the matrix of $n$ transitions $A_n=(a_{ij}(n))$, where
$a_{ij}(n)=A(X_{m+n}=j\vert X_m=i)$
and we get
\begin{equation}\label{chap_kolm_ampl}
a_{ij}(m+n)=\sum_k a_{ik}(m)a_{kj}(n)
\end{equation}
so we have that $A_{m+n}=A_mA_n$ and $A_n=A^n$.

\qee

\section{Probability measures induced by QIFS}\label{s_piqifs}

In this section we present some new results.  Consider a Hilbert space $\mathcal{H}$ of dimension $N=2$. Let $q_1, q_2\in\mathbb{R}$ and also \begin{equation}\label{vi_leva_classico}
V_1=\left(
\begin{array}{cc}
\sqrt{p_{11}} & \sqrt{p_{12}}\\
0 & 0
\end{array}
\right),
\sp V_2=\left(
\begin{array}{cc}
0 & 0\\
\sqrt{p_{21}} & \sqrt{p_{22}}
\end{array}
\right), \sp \rho=\left(
\begin{array}{cc}
\rho_1 & \rho_2\\
\rho_3 & \rho_4
\end{array}
\right)
\end{equation}
We would like to obtain the fixed points for
$$\mathcal{L}(\rho)=q_1V_1\rho V_1^*+q_2V_2\rho V_2^*$$
Then
\begin{equation}\label{sist2mat1}
q_1V_1\rho V_1^*+q_2V_2\rho V_2^*=\rho
\end{equation}
implies
$$q_1\Big[(\sqrt{p_{11}}\rho_1+\sqrt{p_{12}}\rho_3)\sqrt{p_{11}}+(\sqrt{p_{11}}\rho_2+\sqrt{p_{12}}\rho_4)\sqrt{p_{12}}\Big]=\rho_1$$
$$q_2\Big[ (\sqrt{p_{21}}\rho_1+\sqrt{p_{22}}\rho_3)\sqrt{p_{21}}+(\sqrt{p_{21}}\rho_2+\sqrt{p_{22}}\rho_4)\sqrt{p_{22}}\Big]=\rho_4$$
And (\ref{sist2mat1}) also implies that $\rho_2=\rho_3=0$, so we rewrite the system as
$$q_1\Big[\sqrt{p_{11}}\rho_1\sqrt{p_{11}}+\sqrt{p_{12}}\rho_4\sqrt{p_{12}}\Big]=\rho_1$$
$$q_2\Big[ \sqrt{p_{21}}\rho_1\sqrt{p_{21}}+\sqrt{p_{22}}\rho_4\sqrt{p_{22}}\Big]=\rho_4$$
or
\begin{equation}\label{si11}
a\rho_1+f\rho_4=\rho_1
\end{equation}
\begin{equation}\label{si21}
g\rho_1+h\rho_4=\rho_4
\end{equation}
where
$$a=q_1p_{11},\sp f=q_1p_{12},\sp g=q_2p_{21},\sp h=q_2p_{22}$$
We also get that
$$\rho_1=\frac{f}{1-a}\rho_4$$
$$\rho_1=\frac{1-h}{g}\rho_4$$
which is a restriction on the $q_i$, namely
$$\frac{f}{1-a}=\frac{1-h}{g}$$
Therefore the solution of (\ref{si11}) and (\ref{si21}) is
$$\rho=\rho_4\left(
\begin{array}{cc}
\frac{f}{1-a} & 0\\
0 & 1
\end{array}
\right)=\rho_4\left(
\begin{array}{cc}
\frac{1-h}{g} & 0\\
0 & 1
\end{array}
\right)$$ But $\rho_1+\rho_4=1$ implies
\begin{equation}\label{soldc1}
\rho=\left(
\begin{array}{cc}
\frac{q_1p_{12}}{q_1p_{12}-q_1p_{11}+1} & 0\\
0 & \frac{1-q_1p_{11}}{q_1p_{12}-q_1p_{11}+1}
\end{array}
\right)=\left(
\begin{array}{cc}
\frac{1-q_2p_{22}}{1-q_2p_{22}+q_2p_{21}} & 0\\
0 & \frac{q_2p_{21}}{1-q_2p_{22}+q_2p_{21}}
\end{array}
\right)
\end{equation}
Now assume that
$$P=\left(
\begin{array}{cc}
p_{11} & p_{12}\\
p_{21} & p_{22}
\end{array}
\right)$$
is column stochastic. Let $\pi$ be such that $P\pi=\pi$. Such $\pi$ is given by
\begin{equation}\label{solcla1}
\pi=(\frac{p_{12}}{p_{12}-p_{11}+1},\frac{1-p_{11}}{p_{12}-p_{11}+1})
\end{equation}
Compare (\ref{solcla1}) with (\ref{soldc1}). Then fix $q_1=q_2=1$ se we get that
the nonzero entries of $\rho$ are equal to the entries of $\pi$. Such a choice for the $q_i$ is unique. In fact, comparing the $(i,i)$-th entry of $\rho$ with the $i$-th coordinate of $\pi$, we see that if there exists $q_i'$ which make $\rho$ and $\pi$ equal (i.e., the diagonal entries of $\rho$ correspond to the entries of $\pi$), then
$$\frac{q_1p_{12}}{q_1p_{12}-q_1p_{11}+1}=\frac{q_1'p_{12}}{q_1'p_{12}-q_1'p_{11}+1},$$
which implies
$$q_1(q_1'p_{12}-q_1'p_{11}+1)=q_1'(q_1p_{12}-q_1p_{11}+1)$$
$$\Rightarrow q_1q_1'p_{12}-q_1q_1'p_{11}+q_1=q_1q_1'p_{12}-q_1q_1'p_{11}+q_1'$$
and when we cancel terms we get $q_1=q_1'$. In a similar way
$$\frac{1-q_2p_{22}}{1-q_2p_{22}+q_2p_{21}}=\frac{1-q_2'p_{22}}{1-q_2'p_{22}+q_2'p_{21}}$$
implies
$$(1-q_2p_{22})(1-q_2'p_{22}+q_2'p_{21})=(1-q_2'p_{22})(1-q_2p_{22}+q_2p_{21})$$
$$\Rightarrow 1-q_2'p_{22}+q_2'p_{21}-q_2p_{22}+q_2q_2'p_{22}^2-q_2q_2'p_{22}p_{21}$$
$$=1-q_2p_{22}+q_2p_{21}-q_2'p_{22}+q_2q_2'p_{22}^2-q_2q_2'p_{22}p_{21}$$
Then we get
$$q_2'p_{21}=q_2p_{21}\Rightarrow q_2'=q_2$$
and therefore the choice for $q_1$ and $q_2$ is unique.

\qee
Consider a homogeneous QIFS $\mathcal{F}=\{\mathcal{M}_N, F_i,p_i\}_{i=1,\dots, k}$, where
$$F_i(\rho)=\frac{V_i\rho V_i^*}{tr(V_i\rho V_i^*)}$$
 where the $V_i$ are linear with $\sum_i V_i^*V_i=I$ and $p_i(\rho)=tr(V_i\rho V_i^*)$. Then  $\Lambda$ is written as
$$\Lambda(\rho)=\sum_i p_iF_i=\sum_i V_i\rho V_i^*$$
By simplicity we will assume that the quantum system considered can assume two states called $1$ and $2$.

We say that the pair $(\{X_n\}_{n\in\mathbb{N}},\mu)$, $X_n:\Omega\to\{1,\dots, k\}$, is a {\bf Quantum Stochastic Process}, QSP (homogeneous case), associated to the QIFS $\mathcal{F}$ whenever $\mu$ is defined as
\begin{equation}\label{eh_markov_perg}
\mu(X_{1}=x_{1},\dots,X_{n}=x_{n}):=tr(V_{x_{n}}V_{x_{n-1}}\cdots V_{x_{2}}V_{x_{1}}\rho_0 V_{x_{1}}^*V_{x_{2}}^*\cdots V_{x_{n-1}}^*V_{x_{n}}^*)
\end{equation}
where $\rho_0\in\mathcal{M}_N$ is any density operator. The operator $\rho_0$ is a pre-measurement state, that is, we have a quantum system and we prepare $\rho_0$ as being its initial state (for a similar treatment to a sequence of measurements, see the definition of finite dimensional distributions in section \ref{qsp_slom}).

\bigskip

So we can define for any $r$,
\begin{equation}
\mu(X_{r}=x_r\vert X_{{r-1}}=x_{r-1})=\frac{tr(V_{x_r}V_{x_{r-1}}\rho_0 V_{x_{r-1}}^*V_{x_{r}}^*)}{tr(V_{x_{r-1}}\rho_0 V_{x_{r-1}}^*)}
\end{equation}

\begin{defi}
We say that a QSP is {\bf Markov} if
\begin{equation}
\mu(X_n=x_n\vert X_1=x_1,\dots, X_{n-1}=x_{n-1})=\mu(X_n=x_n\vert X_{n-1}=x_{n-1})
\end{equation}
\end{defi}

\qee

{\bf Remark} The condition $\sum_i V_i^* V_i=I$ is enough to show that the measure of a partition of cylinder sets equals 1. For instance, for two states 1 and 2, for $k=2$ and writing
$$\mu(\ov{ij}):=\mu(X_1=i,X_2=j),$$
we have
$$\mu(\ov{11})+\mu(\ov{12})+\mu(\ov{21})+\mu(\ov{22})$$
$$=tr(V_1V_1\rho V_1^* V_1^*)+tr(V_2V_1\rho V_1^* V_2^*)+tr(V_1V_2\rho V_2^* V_1^*)+tr(V_2V_2\rho V_2^* V_2^*)$$
$$=tr(V_1^*V_1[V_1\rho V_1^*])+tr(V_2^*V_2[V_1\rho V_1^*] )+tr(V_1^*V_1[V_2\rho V_2^*])+tr(V_2^*V_2[V_2\rho V_2^*])$$
$$=tr\Big((V_1^*V_1+V_2^*V_2)[V_1\rho V_1^*]\Big)+tr\Big((V_1^*V_1+V_2^*V_2)[V_2\rho V_2^*]\Big)$$
\begin{equation}
=tr(V_1\rho V_1^*)+tr(V_2\rho V_2^*)=tr((V_1^*V_1+V_2^*V_2)\rho)=1
\end{equation}
However, we note that there exist examples in which we can show that the measure of a partition of cylinder sets equals 1 even if we do not suppose that $\sum_i V_i^* V_i=I$. This happens, for instance, in the following construction involving stochastic matrices.

% ver 010809.mws

\qee

Let us consider the particular case in which the operator $\rho_0\in\mathcal{M}_N$, given in the definition of QSP is a fixed point for $\Lambda(\rho)=\sum_{i=1}^k V_i\rho V_i^*$ induced by the QIFS $\mathcal{F}$.

Suppose that $V_1$ and $V_2$ are defined by (\ref{vi_leva_classico}). Suppose that the matrix $P=(p_{ij})$ is column stochastic and that we have $\pi$ such that $P\pi=\pi$. For instance we have
\begin{equation}
\mu(X_1=1,X_2=2)=tr(V_2 V_1\rho_0 V_1^*V_2^*)=p_{21}(p_{11}\rho_{11}+p_{12}\rho_{22})=p_{21}\rho_{11}
\end{equation}
because with the choice of $V_i$ we made, we have that the nonzero entries of $\rho_0$ correspond to the entries of $\pi$. So we can interpret $p_{ij}$ as being
\begin{equation}
p_{ij}=\mu(X_2=j\vert X_1=i)
\end{equation}
In a similar way,
\begin{equation}
\mu(X_1=2,X_2=1)=tr(V_1 V_2\rho_0 V_2^*V_1^*)=p_{12}\rho_{22}
\end{equation}
and
\begin{equation}
\mu(X_1=1,X_2=2,X_3=1)=tr(V_1 V_2 V_1\rho_0 V_1^*V_2^*V_1^*)=p_{12}p_{21}\rho_{11}
\end{equation}

\bigskip

{\bf Remark} A simple calculation shows that with the $V_i$ given by (\ref{vi_leva_classico}) we have that $\sum_i V_i^* V_i\neq I$. However, we still have that
$$\mu(\ov{11})+\mu(\ov{12})+\mu(\ov{21})+\mu(\ov{22})=1$$

\qee

To prove that the choice (\ref{vi_leva_classico}) reduces to the classic case for any sequence, we use the following lemma.

\begin{lem}\label{lema_comomult_vi}
Suppose $N=2$, $k=2$. Then for every $m$, for $V_i$ given by (\ref{vi_leva_classico}) and $\rho_0$ corresponding to the stationary vector $\pi$ for $P$, we have that the product
\begin{equation}
V_{x_m}V_{x_{m-1}}\cdots V_{x_1}\rho_0 V_{x_1}^*V_{x_2}^*\cdots V_{x_m}^*
\end{equation}
has the form
\begin{equation}
\left(
\begin{array}{cc}
* & 0\\
0 & 0
\end{array}
\right) \sp\textrm{ ou }\sp \left(
\begin{array}{cc}
0 & 0\\
0 & *
\end{array}
\right)
\end{equation}
depending on whether $x_m=1$ or $x_m=2$, respectively.
\end{lem}
{\bf Proof} By induction. If $m=1$ then
\begin{equation}
V_1\rho_0 V_1^*=\left(
\begin{array}{cc}
p_{11}\rho_{11}+p_{12}\rho_{22} & 0\\
0 & 0
\end{array}
\right)
\end{equation}
and
\begin{equation}
V_2\rho_0 V_2^*=\left(
\begin{array}{cc}
0 & 0\\
0 & p_{21}\rho_{11}+p_{22}\rho_{22}
\end{array}
\right)
\end{equation}
Suppose the lemma valid for $m$, we consider the product
\begin{equation}
V_{x_{m+1}}V_{x_{m}}\cdots V_{x_1}\rho_0 V_{x_1}^*V_{x_2}^*\cdots V_{x_m}^*V_{x_{m+1}}^*
\end{equation}
Suppose $x_{m+1}=1$. Then a simple calculation shows that
\begin{equation}
V_{1}\left(
\begin{array}{cc}
* & 0\\
0 & 0
\end{array}
\right) V_{1}^*\sp\textrm{ e }\sp V_{1}\left(
\begin{array}{cc}
0 & 0\\
0 & *
\end{array}
\right) V_{1}^*
\end{equation}
has only one nonzero entry, namely the $(1,1)$-th entry. We proceed in a similar way for the case that $x_{m+1}=2$, that is
\begin{equation}
V_{2}\left(
\begin{array}{cc}
* & 0\\
0 & 0
\end{array}
\right) V_{2}^*\sp\textrm{ e }\sp V_{2}\left(
\begin{array}{cc}
0 & 0\\
0 & *
\end{array}
\right) V_{2}^*
\end{equation}
has only one nonzero entry, namely the $(2,2)$-th entry.

\qed

\begin{pro}
If we set
\begin{equation}\label{vi_leva_classico2}
V_1=\left(
\begin{array}{cc}
\sqrt{p_{11}} & \sqrt{p_{12}}\\
0 & 0
\end{array}
\right),
\sp V_2=\left(
\begin{array}{cc}
0 & 0\\
\sqrt{p_{21}} & \sqrt{p_{22}}
\end{array}
\right),
\end{equation}
then
\begin{equation}
\mu(X_1=x_1,X_2=x_2,\dots ,X_n=x_n)=p_{x_{n}x_{n-1}}p_{x_{n-1}x_{n-2}}\cdots p_{x_3x_2}p_{x_2x_1}\rho_{x_1x_1}
\end{equation}
where $\rho_{ij}$ denotes the $(i,j)$-th entry of $\rho_0$, eigenstate for $\Lambda(\rho)=\sum_i V_i\rho V_i^*$.
\end{pro}
{\bf Proof} We prove by induction. Suppose $n=1$. Then
$$\mu(X_{1}=1)=tr(V_1\rho_0 V_1^*)=p_{11}\rho_{11}+p_{12}\rho_{22}=\rho_{11}$$
$$\mu(X_{1}=2)=tr(V_2\rho_0 V_2^*)=p_{21}\rho_{11}+p_{22}\rho_{22}=\rho_{22}$$
For the sake of clarity we also show the case $n=2$. We have, after some routine calculations that
\begin{equation}
\mu(X_{1}=1,X_{2}=1)=tr(V_1V_1\rho_0 V_1^*V_1^*)=p_{11}\rho_{11}
\end{equation}
\begin{equation}
\mu(X_{1}=1,X_{2}=2)=tr(V_2V_1\rho_0 V_1^*V_2^*)=p_{21}\rho_{11}
\end{equation}
\begin{equation}
\mu(X_{1}=2,X_{2}=1)=tr(V_1V_2\rho_0 V_2^*V_1^*)=p_{12}\rho_{22}
\end{equation}
\begin{equation}
\mu(X_{1}=2,X_{2}=2)=tr(V_2V_2\rho_0 V_2^*V_2^*)=p_{22}\rho_{22}
\end{equation}
Now suppose the lemma holds for $n$, let us prove it for $n+1$.

\bigskip

First suppose $x_{n+1}=1$. Then
$$\mu(X_{1}=x_{1},\dots,X_{n}=x_{n},X_{{n+1}}=1)$$
\begin{equation}
=tr(V_{1}V_{x_{n}}\cdots V_{x_{2}}V_{x_{1}}\rho_0 V_{x_{1}}^*V_{x_{2}}^*\cdots V_{x_{n-1}}^*V_{x_{n}}^*V_{1}^*)
\end{equation}
Using lemma \ref{lema_comomult_vi}, we have two cases. If $x_n=1$ then
$$V_{x_{n}}\cdots V_{x_{2}}V_{x_{1}}\rho_0 V_{x_{1}}^*V_{x_{2}}^*\cdots V_{x_{n-1}}^*V_{x_{n}}^*=\left(
\begin{array}{cc}
* & 0\\
0 & 0
\end{array}
\right)$$
and therefore
$$V_{1}V_{x_{n}}\cdots V_{x_{2}}V_{x_{1}}\rho_0 V_{x_{1}}^*V_{x_{2}}^*\cdots V_{x_{n-1}}^*V_{x_{n}}^*V_{1}^*=V_1\left(
\begin{array}{cc}
* & 0\\
0 & 0
\end{array}
\right) V_1^*=\left(
\begin{array}{cc}
*p_{11} & 0\\
0 & 0
\end{array}
\right)$$
and so by taking the trace we get
$$tr(V_{1}V_{x_{n}}\cdots V_{x_{2}}V_{x_{1}}\rho_0 V_{x_{1}}^*V_{x_{2}}^*\cdots V_{x_{n-1}}^*V_{x_{n}}^*V_{1}^*)$$
\begin{equation}
= p_{11}p_{1x_{n-1}}p_{x_{n-1}x_{n-2}}\cdots p_{x_3x_2}p_{x_2x_1}\rho_{x_1x_1}
\end{equation}
In a similar way, if $x_n=2$,
$$V_{x_{n}}\cdots V_{x_{2}}V_{x_{1}}\rho_0 V_{x_{1}}^*V_{x_{2}}^*\cdots V_{x_{n-1}}^*V_{x_{n}}^*=\left(
\begin{array}{cc}
0 & 0\\
0 & *
\end{array}
\right)$$
$$V_{1}V_{x_{n}}\cdots V_{x_{2}}V_{x_{1}}\rho_0 V_{x_{1}}^*V_{x_{2}}^*\cdots V_{x_{n-1}}^*V_{x_{n}}^*V_{1}^*=V_1\left(
\begin{array}{cc}
0 & 0\\
0 & *
\end{array}
\right) V_1^*=\left(
\begin{array}{cc}
*p_{12} & 0\\
0 & 0
\end{array}
\right)$$
and taking the trace gives
$$tr(V_{1}V_{x_{n}}\cdots V_{x_{2}}V_{x_{1}}\rho_0 V_{x_{1}}^*V_{x_{2}}^*\cdots V_{x_{n-1}}^*V_{x_{n}}^*V_{1}^*)$$
\begin{equation}
= p_{12}p_{2x_{n-1}}p_{x_{n-1}x_{n-2}}\cdots p_{x_3x_2}p_{x_2x_1}\rho_{x_1x_1}
\end{equation}

\bigskip

Now we suppose $x_{n+1}=2$, and we proceed in an analogous way.
$$\mu(X_{1}=x_{1},\dots,X_{n}=x_{n},X_{{n+1}}=2)$$
\begin{equation}
=tr(V_{2}V_{x_{n}}\cdots V_{x_{2}}V_{x_{1}}\rho_0 V_{x_{1}}^*V_{x_{2}}^*\cdots V_{x_{n-1}}^*V_{x_{n}}^*V_{2}^*)
\end{equation}
By lemma \ref{lema_comomult_vi}, we have two cases. If $x_n=1$ then
$$V_{x_{n}}\cdots V_{x_{2}}V_{x_{1}}\rho_0 V_{x_{1}}^*V_{x_{2}}^*\cdots V_{x_{n-1}}^*V_{x_{n}}^*=\left(
\begin{array}{cc}
* & 0\\
0 & 0
\end{array}
\right)$$
therefore
$$V_{2}V_{x_{n}}\cdots V_{x_{2}}V_{x_{1}}\rho_0 V_{x_{1}}^*V_{x_{2}}^*\cdots V_{x_{n-1}}^*V_{x_{n}}^*V_{2}^*=V_2\left(
\begin{array}{cc}
* & 0\\
0 & 0
\end{array}
\right) V_2^*=\left(
\begin{array}{cc}
*p_{21} & 0\\
0 & 0
\end{array}
\right)$$
and taking the trace we get
$$tr(V_{2}V_{x_{n}}\cdots V_{x_{2}}V_{x_{1}}\rho_0 V_{x_{1}}^*V_{x_{2}}^*\cdots V_{x_{n-1}}^*V_{x_{n}}^*V_{2}^*)$$
\begin{equation}
= p_{21}p_{1x_{n-1}} p_{x_{n-1}x_{n-2}}\cdots p_{x_3x_2}p_{x_2x_1}\rho_{x_1x_1}
\end{equation}
Analogously if $x_n=2$
$$V_{x_{n}}\cdots V_{x_{2}}V_{x_{1}}\rho_0 V_{x_{1}}^*V_{x_{2}}^*\cdots V_{x_{n-1}}^*V_{x_{n}}^*=\left(
\begin{array}{cc}
0 & 0\\
0 & *
\end{array}
\right)$$
$$V_{2}V_{x_{n}}\cdots V_{x_{2}}V_{x_{1}}\rho_0 V_{x_{1}}^*V_{x_{2}}^*\cdots V_{x_{n-1}}^*V_{x_{n}}^*V_{2}^*=V_2\left(
\begin{array}{cc}
0 & 0\\
0 & *
\end{array}
\right) V_2^*=\left(
\begin{array}{cc}
*p_{22} & 0\\
0 & 0
\end{array}
\right)$$
and taking the trace
$$tr(V_{2}V_{x_{n}}\cdots V_{x_{2}}V_{x_{1}}\rho_0 V_{x_{1}}^*V_{x_{2}}^*\cdots V_{x_{n-1}}^*V_{x_{n}}^*V_{2}^*)$$
\begin{equation}
= p_{22} p_{2x_{n-1}}p_{x_{n-1}x_{n-2}}\cdots p_{x_3x_2}p_{x_2x_1}\rho_{x_1x_1}
\end{equation}
\qed

\begin{cor}
The quantum stochastic process induced by
\begin{equation}
V_1=\left(
\begin{array}{cc}
\sqrt{p_{11}} & \sqrt{p_{12}}\\
0 & 0
\end{array}
\right),
\sp V_2=\left(
\begin{array}{cc}
0 & 0\\
\sqrt{p_{21}} & \sqrt{p_{22}}
\end{array}
\right),
\end{equation}
is Markov.
\end{cor}
{\bf Proof} By the proposition, we have that the measure $\mu$ reduces to the Markov measure for matrices.

\qed

\begin{lem}\label{lema_stationary_time}
For $V_i$ linear maps and $\rho_0$ fixed point for $\Lambda=\sum_i V_i\rho V_i^*$, we have for any $m, n$,
$$\mu(X_1=x_1, X_2=x_2, \dots ,X_n=x_n)=\mu(X_m=x_1, X_{m+1}=x_2, \dots ,X_{m+n}=x_n)$$
\end{lem}
{\bf Proof} We prove the lemma for the case in which we have two possible states $1$ and $2$. We have
$$\mu(X_m=x_1, X_{m+1}=x_2, \dots ,X_{m+n}=x_n)$$
$$=\sum_{i_1,\dots ,i_{m-1}}\mu(X_{1}=i_1,X_2=i_2,\dots, X_{m-1}=i_{m-1},X_m=x_1,\dots ,X_{m+n}=x_n)$$
$$=\sum_{i_2,\dots ,i_{m-1}}tr(V_{x_n}\cdots V_{x_1}V_{i_{m-1}}\cdots V_{i_2}V_1\rho_0 V_1^* V_{i_2}^*\cdots ) $$
$$+ tr(V_{x_n}\cdots V_{x_1}V_{i_{m-1}}\cdots V_{i_2}V_2\rho_0 V_2^* V_{i_2}^*\cdots )$$
$$=\sum_{i_2,\dots ,i_{m-1}}tr(V_{x_n}\cdots V_{x_1}V_{i_{m-1}}\cdots V_{i_2}\rho_0 V_{i_2}^*V_{i_3}^*\cdots  V_{i_{m-1}}^*V_{x_1}^*\cdots V_{x_n}^*)$$
Repeating the procedure above for $i_2$, $i_3$, etc. we get
$$\mu(X_m=x_1, X_{m+1}=x_2, \dots ,X_{m+n}=x_n)=tr(V_{x_n}\cdots V_{x_1}\rho_0 V_{x_1}^*\cdots V_{x_n}^*)$$
This concludes the proof.
\qed

\begin{examp}
Let us make an inspection with respect to the Chapman-Kolmogorov equation, that is, we would like to know if the equality
\begin{equation}
\mu_{ij}(m+n)=\sum_k \mu_{ik}(m)\mu_{kj}(n)
\end{equation}
holds, where
$$\mu_{ij}(n)=\mu(X_{m+n}=j\vert X_m=i)$$
Take for instance, $m=n=i=j=1$. Then
$$\sum_k \mu_{ik}(m)\mu_{kj}(n)=\mu_{11}(1)\mu_{11}(1)+ \mu_{12}(1)\mu_{21}(1)$$
\begin{equation}\label{no_ck1}
=\frac{tr(V_1V_1\rho V_1^*V_1^*)^2}{tr(V_1\rho V_1^*)^2}+\frac{tr(V_2V_1\rho V_1^*V_2^*)}{tr(V_1\rho V_1^*)}\frac{tr(V_1V_2\rho V_2^*V_1^*)}{tr(V_2\rho V_2^*)}
\end{equation}
and
$$\mu_{ij}(m+n)=\mu_{11}(2)=\mu(X_3=1\vert X_1=1)$$
\begin{equation}\label{no_ck2}
=\frac{tr(V_1V_1V_1\rho V_1^*V_1^*V_1^*)}{tr(V_1\rho V_1^*)}+\frac{tr(V_1V_2V_1\rho V_1^*V_2^*V_1^*)}{tr(V_1\rho V_1^*)}
\end{equation}
Now let $V_1$, $V_2$ be given by (\ref{vi_leva_classico2}), then we obtain classic calculations, so the Chapman-Kolmogorov equation holds. Now take
\begin{equation}\label{um_nao_markov}
V_1=\left(
\begin{array}{cc}
1 & 0\\
0 & 0
\end{array}
\right),\sp V_2=\left(
\begin{array}{cc}
1 & 0\\
0 & 2
\end{array}
\right)
\end{equation}
then we get, from (\ref{no_ck1}) and (\ref{no_ck2}):
\begin{equation}
\frac{tr(V_1V_1\rho V_1^*V_1^*)^2}{tr(V_1\rho V_1^*)^2}+\frac{tr(V_2V_1\rho V_1^*V_2^*)}{tr(V_1\rho V_1^*)}\frac{tr(V_1V_2\rho V_2^*V_1^*)}{tr(V_2\rho V_2^*)}=1+\frac{\rho_{11}}{\rho_{11}+4\rho_{22}}
\end{equation}
and
\begin{equation}
\frac{tr(V_1V_1V_1\rho V_1^*V_1^*V_1^*)}{tr(V_1\rho V_1^*)}+\frac{tr(V_1V_2V_1\rho V_1^*V_2^*V_1^*)}{tr(V_1\rho V_1^*)}=1+1=2
\end{equation}
% ver calculos em 280709.mws, 210909 also
Then in this case we have that the Chapman-Kolmogorov equation holds if and only if $\rho_{22}=0$ that is, if $\rho_{11}=1$. Also, we note that $\sum_i V_i^* V_i\neq I$.
To conclude this example, we take $V_1$ and $V_2$ with $\sum_i V_i^*V_i=I$, namely,
\begin{equation}
V_1=\left(
\begin{array}{cc}
\frac{1}{\sqrt{3}} & 0\\
0 & 0
\end{array}
\right),\sp V_2=\left(
\begin{array}{cc}
\sqrt{\frac{2}{3}} & 0\\
0 & 1
\end{array}
\right)
\end{equation}
Take for instance $\rho_0=\frac{1}{4}\vert 1\rangle\langle 1\vert+\frac{3}{4}\vert 2\rangle\langle 2\vert$, a fixed point for the associated $\Lambda$. A simple calculation shows that (\ref{no_ck1}) and (\ref{no_ck2}) are different. Therefore our calculation shows that the Chapman-Kolmogorov equation does not hold in general (for our setting).

\end{examp}

\qee

We would like to obtain a nonhomogeneous version for the measure we defined in
(\ref{eh_markov_perg}) in the homogeneous case, i.e., we are looking for a measure induced by a nonhomogeneous QIFS. Let $W_i$, $i=1,\dots ,k$ be linear operators such that $\sum_i W_i^*W_i=I$. Let $\rho_0\in\mathcal{M}_N$. Define
$$\mu(X_{1}=x_{1},\dots,X_{n}=x_{n}):=$$
$$=tr(W_{x_1}\rho_0 W_{x_1}^*)\frac{tr(W_{x_2}V_{x_1}\rho_0 V_{x_1}^*W_{x_2}^*)}{tr(V_{x_1}\rho_0V_{x_1}^*)}\frac{tr(W_{x_3}V_{x_2}V_{x_1}\rho_0 V_{x_1}^*V_{x_2}^*W_{x_3}^*)}{tr(V_{x_2}V_{x_1}\rho_0V_{x_1}^*V_{x_2}^*)}\times\cdots$$
$$\cdots\times\frac{tr(W_{x_{n-1}}V_{x_{n-2}}\cdots V_{x_1}\rho_0 V_{x_1}^*\cdots V_{x_{n-2}}^*W_{x_{n-1}}^*)}{tr(V_{x_{n-2}}\cdots V_{x_1}\rho_0V_{x_1}^*\cdots V_{x_{n-2}}^*)}\times$$
\begin{equation}
\times\frac{tr(W_{x_{n}}V_{x_{n-1}}\cdots V_{x_1}\rho_0 V_{x_1}^*\cdots V_{x_{n-1}}^*W_{x_{n}}^*)}{tr(V_{x_{n-1}}\cdots V_{x_1}\rho_0V_{x_1}^*\cdots V_{x_{n-1}}^*)}
\end{equation}
that is,
$$\mu(X_{1}=x_{1},\dots,X_{n}=x_{n}):=$$
\begin{equation}\label{escolha_poss_1}
tr(W_{x_1}\rho_0 W_{x_1}^*)\prod_{i=2}^n \frac{tr(W_{x_i}V_{x_{i-1}}\cdots V_{x_1}\rho_0 V_{x_1}^*\cdots V_{x_{i-1}}^*W_{x_i}^*)}{tr(V_{x_{i-1}}V_{x_{i-2}}\cdots V_{x_1}\rho_0 V_{x_1}^*\cdots V_{x_{i-2}}^*V_{x_{i-1}}^*)}
\end{equation}

\bigskip

{\bf Remark} A calculation shows that if we suppose $\sum_i W_i^* W_i=I$, then
$$\sum_{i_1,\dots i_n} \mu(\ov{i_1\cdots i_n})=1$$
Besides, if we suppose that $W_i=V_i$ for all $i$, then we recover the measure definition for homogeneous QSP.

\qee

Consider a QIFS $\mathcal{F}=\{\mathcal{M}_N, F_i,p_i\}_{i=1,\dots, k}$, where
$$F_i(\rho)=\frac{V_i\rho V_i^*}{tr(V_i\rho V_i^*)}$$
where the $V_i$ are linear and $p_i(\rho)=tr(W_i\rho W_i^*)$, com $\sum_i W_i^*W_i=I$.

\bigskip

\begin{defi}
We say that the pair $(\{X_n\}_{n\in\mathbb{N}},\mu)$, $X_n:\Omega\to\{1,\dots, k\}$, is a {\bf Quantum Stochastic Process} associated to the {\bf nonhomogeneous} QIFS $\mathcal{F}$ if $\mu$ is defined by (\ref{escolha_poss_1}), where $\rho_0\in\mathcal{M}_N$ is any density operator.
\end{defi}

{\bf Remark} In the definition above we can, of course, consider the particular case in which   $\rho_0$ is a fixed point for
$$\Lambda(\rho)=\sum_{i=1}^k tr(W_i\rho W_i^*)\frac{V_i\rho V_i^*}{tr(V_i\rho V_i^*)},$$
induced by the QIFS $\mathcal{F}$.

\qee

Recall that by lemma \ref{lema_stationary_time}, a homogeneous QSP is always stationary. This is no longer true in general for nonhomogeneous QSP.

\begin{examp}
Let $\{X_n\}_{n\in\mathbb{N}}$ be a QSP induced by a nonhomogeneous QIFS. We would like to know whether
\begin{equation}
\mu(X_1=1,X_2=2)=\mu(X_2=1,X_3=2)
\end{equation}
By definition we have:
\begin{equation}
\mu(X_1=1,X_2=2)=tr(W_1\rho_0 W_1^*)\frac{tr(W_2 V_1\rho_0 V_1^* W_2^*)}{tr(V_1\rho_0 V_1^*)}
\end{equation}
And also
$$\mu(X_2=1,X_3=2)=\mu(X_1=1,X_2=1,X_3=2)+\mu(X_1=2,X_2=1,X_3=2)$$
$$=tr(W_1\rho_0 W_1^*)\frac{tr(W_1 V_1\rho_0 V_1^*W_1^*)}{tr(V_1\rho_0 V_1^*)}\frac{tr(W_2V_1 V_1\rho_0 V_1^*V_1^*W_2^*)}{tr(V_1V_1\rho_0 V_1^*V_1^*)}$$
\begin{equation}
+ tr(W_2\rho_0 W_2^*)\frac{tr(W_1 V_2\rho_0 V_2^*W_1^*)}{tr(V_2\rho_0 V_2^*)}\frac{tr(W_2V_1 V_2\rho_0 V_2^*V_1^*W_2^*)}{tr(V_1V_2\rho_0 V_2^*V_1^*)}
\end{equation}

\bigskip

$$=tr\Bigg[W_2 V_1\Big[tr(W_1\rho_0 W_1^*)\frac{V_1\rho_0 V_1^*}{tr(V_1\rho_0 V_1^*)}\Big(\frac{tr(W_1 V_1\rho_0 V_1^*W_1^*)}{tr(V_1V_1\rho_0 V_1^*V_1^*)}\Big) +$$
\begin{equation}\label{ep_nest1}
+tr(W_2\rho_0 W_2^*)\frac{V_2\rho_0 V_2^*}{tr(V_2\rho_0 V_2^*)}\Big(\frac{tr(W_1 V_2\rho_0 V_2^*W_1^*)}{tr(V_1V_2\rho_0 V_2^*V_1^*)}\Big)\Big]V_1^*W_2^* \Bigg]
\end{equation}
Note that in the homogeneous case we have that both fractions in parenthesis on equation (\ref{ep_nest1}) are equal to 1, so if $\rho_0$ is a fixed point for $\Lambda$, then we have stationarity, a fact we have already proved. But in the nonhomogeneous case, the terms in parenthesis are not equal to 1 in general.

\end{examp}

\qee

\section{A definition of entropy for QIFS}\label{entr_artur}

We will present a notion of entropy for ``invariant" (or
``stationary") measures with support on density matrices. This
definition is obtained by adapting the reasoning described in
\cite{GL}, \cite{lopes} and \cite{lopes_elismar}  to the present
situation. The main idea is to define this concept via the Ruelle
operator and to avoid the use of partitions.

Denote by $p$ an arbitrary choice of mappings $p_i:\mathcal{M}_N\to\mathbb{R}$,
$i=1,\dots, k$ for a certain $k$. Let
$$m_b(\mathcal{M}_N):=\{f:\mathcal{M}_N\to\mathbb{R} : \textrm{f is measurable and bounded}\}$$
Let $\mathcal{U}_p:m_b(\mathcal{M}_N)\to m_b(\mathcal{M}_N)$,
$$(\mathcal{U}_pf)(\rho):=\sum_{i=1}^k p_i(\rho)f(F_i(\rho))$$
Let us consider all possible choices of mappings $p_i:\mathcal{M}_N\to\mathbb{R}$  which satisfy
\begin{equation}\label{pb_pap1}
\mathcal{U}_p \,1\,=1
\end{equation}
Each $p$ determines an operator $\mathcal{U}_p$. The set of all possible $p$ that satisfy (\ref{pb_pap1}) will be denoted by $P$.

Let $(\mathcal{M}_N,F_i,p_i)_{i=1,\dots k}$ be a QIFS. An example of Markov operator for measures is the one we defined before, given by $\mathcal{V}_p : M^1(\mathcal{M}_N)\to M^1(\mathcal{M}_N)$,
$$(\mathcal{V}_p\nu)(B)=\sum_{i=1}^k\int_{F_i^{-1}(B)}p_id\nu,$$
which we will call {\bf the Markov operator Markov induced by the} $p_i$. That is, we will consider all $\mathcal{V}_p$ with $p\in P$. We say that $\nu$ is {\bf invariant} for the $F_i$ if for some $p\in P$ we have that $\mathcal{V}_p\nu=\nu$.

\bigskip

Let ${\mathcal M}_F$ be the set of all invariant measures for a
fixed choice of the dynamics $F_i$, $i=1,\dots, k$. For such
measures $\nu\in{\mathcal M}_F$, and based on \cite{GL}, \cite{lopes} and
\cite{lopes_elismar}, define
$$h_0(\nu):=\inf_{f \in \mathbb{B}^{+}}\int \log(\sum_{i=1}^{k} \frac{ f \circ F_{i}}{f}) d\nu$$
Above, $\mathbb{B}^{+}$ denotes the bounded, positive, borelean functions on $\mathcal{M}_N$.

\begin{pro}\label{entropia>0}
For $\nu\in\mathcal{M}_F$, we have that $0\leq h_0(\nu)\leq\log{k}$.
\end{pro}

In order to prove this proposition, we need the following lemma.

\begin{lem}\label{lema_log_conv}\cite{lopes_elismar}
Let $\beta\geq 1+\alpha$ and numbers $a_i\in [1+\alpha,\beta]$, $i=1,\dots, k$. Then there exists $\epsilon\geq 1$ such that
$$\log{\Big(\epsilon\sum_{i=1}^k a_i\Big)}\geq \sum_{i=1}^k \log{(\epsilon a_i)}.$$
\end{lem}

The proof of this lemma follows by choosing
$$\epsilon=\exp{\Big(\frac{1}{k}\frac{\log{\sum_{i=1}^k a_i}}{\sum_{i=1}^k\log{a_i}} \Big)}$$

\begin{lem}\label{lema_sifcn}
If $f\in\mathbb{B}^+$ and $\nu\in\mathcal{M}_F$ then
$$\sum_{i=1}^k\int f\circ F_i d\nu\geq \int f d\nu$$
\end{lem}
{\bf Proof}
First suppose that $f=1_B$, where $B$ is a measurable set. We have that
$$\sum_{i=1}^k\int 1_B\circ F_i d\nu\geq\sum_{i=1}^k\int p_i(x)1_B(F_i(x))d\nu(x)=\sum_{i=1}^k\int_{F_i^{-1}(B)}p_i(x)d\nu(x)$$
$$=\mathcal{V}_p(\nu)(B)=\nu(B)=\int 1_Bd\nu$$
Then, assume that $f=\sum_{j=1}^l b_j 1_{B_j}$, i.e., a simple function. Then
$$\sum_{i=1}^k\int\sum_{j=1}^l b_j1_{B_j}\circ F_id\nu=\sum_{j=1}^l b_j\sum_{i=1}^k\int 1_{B_j}\circ F_id\nu$$
$$\geq \sum_{j=1}^l b_j\sum_{i=1}^k\int p_i(x)1_{B_j}(F_i(x)) d\nu=\sum_{j=1}^l b_j \mathcal{V}_p(\nu)(B_j)$$
$$=\sum_{j=1}^l b_j\nu(B_j)=\int fd\nu$$
Now let $f=\lim_n f_n$, a limit of a sequence of simple functions. Note that we suppose $f\in\mathbb{B}^+$, so $f$ is bounded, and since $\nu$ is a probability measure on $\mathcal{M}_N$, it follows that $f$ is integrable. By the bounded convergence theorem, we have that
$$\sum_{i=1}^k\int f\circ F_id\nu=\sum_{i=1}^k\int \lim_n{f_n\circ F_i}d\nu=\lim_n \sum_{i=1}^k\int f_n\circ F_id\nu$$
$$\geq \lim_n \int f_nd\nu=\int\lim_n f_nd\nu=\int fd\nu$$

\qed

The following proof is an adaptation of results seen in \cite{lopes_elismar}.

\vspace{0.2cm}

{\bf Proof of proposition \ref{entropia>0}} Let us restrict the proof  for the case in which we have a QIFS $(\mathcal{M}_N, F_i, p_i)_{i=1,\dots, k}$, where $F_i(\rho)=V_i\rho V_i^*$, with linear $V_i$.

\vspace{0.2cm}

First note that if $f\equiv 1$, we have $\int\log(\sum_{i=1}^k 1)d\nu=\log{k}$, so $h_0(\nu)\leq\log{k}$.

Let $I=\int\log{(\sum_{i=1}^k\frac{f\circ F_i}{f})}d\nu$ and
suppose, without loss of generality, that $1+\alpha\leq f\leq \beta$
(note that this integral is invariant by the projective mapping $f\to \lambda f$). Then
\begin{equation}\label{eli_eq1}
I=\int\log{(\sum_{i=1}^k\frac{\epsilon f\circ F_i}{\epsilon f})}d\nu=\int\log{(\sum_{i=1}^k\epsilon f\circ F_i)}d\nu-\int\log{(\epsilon f)}d\nu
\end{equation}
Define
$$a_i=f\circ F_i(\rho)$$
Then
$$\epsilon(\rho)=\exp{\Big(\frac{1}{k}\frac{\log{\sum_{i=1}^k f\circ F_i}}{\sum_{i=1}^k\log{f\circ F_i}} \Big)}\geq\epsilon_0\geq 1,$$
by the compactness of $\mathcal{M}_N$. With such choice we obtain, by lemma (\ref{lema_log_conv}),
\begin{equation}\label{eli_eq2}
\log{(\epsilon_0\sum_{i=1}^k f\circ F_i)}\geq \sum_{i=1}^k\log{(\epsilon_0 f\circ F_i)}
\end{equation}
Apply (\ref{eli_eq2}) on (\ref{eli_eq1}), then
$$I\geq \sum_{i=1}^k\int\log{(\epsilon_0 f\circ F_i)}d\nu-\int\log{(\epsilon_0 f)}d\nu$$
Then by lemma (\ref{lema_sifcn}) applied on the function $\log{(\epsilon f)}$ (note that we have $\log{(\epsilon_0 f)}\in\mathbb{B}^+$ because $\epsilon_0\geq 1)$, we get
$$I\geq\int\log{(\epsilon f)}d\nu - \int\log{(\epsilon f)}d\nu=0$$
\qed

The computation in the next example shows that the concept of entropy described here is different from the one presented in \cite{BLLT1} \cite{BLLT2}.

\begin{examp} We will consider an example of a probability $\eta$ such that
$\mathcal{V}(\eta)
=\eta$ and we will compute the entropy of $\eta$.

Suppose a QIFS, such that
$$p_i(\rho)=tr(W_i\rho W_i^*),\sp \sum_i W_i^* W_i=I,\sp F_i(\rho)=\frac{V_i\rho V_i^*}{tr(V_i\rho V_i^*)}$$
 for $i=1,\dots, k$. Denote $m_b(\mathcal{M}_N)$ the space of bounded and measurable functions in $\mathcal{M}_N$.

 Consider  $\Lambda:\mathcal{M}_N\to \mathcal{M}_N$,
$$\Lambda(\rho)=\sum_i p_i(\rho)F_i(\rho)=\sum_i tr(W_i\rho W_i^*)\frac{V_i\rho V_i^*}{tr(V_i\rho V_i^*)}$$

 Suppose there exists  a density matrix $\rho$ which $\Lambda$-invariant. As we know, such state is the barycenter of  $\mu$ which is  $\mathcal{V}$-invariant \cite{BLLT1}.

Suppose $\mathcal{V}\mu=\mu$, then we can write
$$\int f d\mu=\int f d\mathcal{V}\mu=\sum_{i=1}^k \int p_i(\rho)f(F_i(\rho))d\mu(\rho)=\sum_i\int p_i(\rho)f\Big(\frac{V_i\rho V_i^*}{tr(V_i\rho V_i^*)}\Big)d\mu$$
$$=\sum_i\int tr(W_i\rho W_i^*)f\Big(\frac{V_i\rho V_i^*}{tr(V_i\rho V_i^*)}\Big)d\mu$$

Therefore, for any $f\in m_b(\mathcal{M}_N)$, we got the condition
\begin{equation}\label{n_eq_set09}
\int f d\mu=\sum_i\int tr(W_i\rho W_i^*)f\Big(\frac{V_i\rho V_i^*}{tr(V_i\rho V_i^*)}\Big)d\mu
\end{equation}

Let us consider a particular example where $N=2$, $k=4$, and
$$V_1=\left(
\begin{array}{cc}
\sqrt{p_{11}} & 0\\
0 & 0
\end{array}
\right),
\sp V_2=\left(
\begin{array}{cc}
0 & \sqrt{p_{12}}\\
0 & 0
\end{array}
\right)
,$$
$$V_3=\left(
\begin{array}{cc}
0 & 0\\
\sqrt{p_{21}} & 0
\end{array}
\right),
\sp V_4=\left(
\begin{array}{cc}
0 & 0\\
0 & \sqrt{p_{22}}
\end{array}
\right),
$$
in such way that the $p_{ij}$ are the entries of a column stochastic  matrix $P$. Let $\pi=(\pi_1,\pi_2)$ be a vector such that  $P\pi=\pi$. A simple calculation shows that for $\rho$, the density matrix such that has entries $\rho_{ij}$, we have
\begin{equation}\label{300909a}
V_1\rho V_1^*=\left(
\begin{array}{cc}
p_{11}\rho_{11} & 0\\
0 & 0
\end{array}
\right),\sp V_2\rho V_2^*=\left(
\begin{array}{cc}
p_{12}\rho_{22} & 0\\
0 & 0
\end{array}
\right)
\end{equation}
\begin{equation}\label{300909b}
V_3\rho V_3^*=\left(
\begin{array}{cc}
0 & 0\\
0 & p_{21}\rho_{11}
\end{array}
\right),\sp V_4\rho V_4^*=\left(
\begin{array}{cc}
0 & 0\\
0 & p_{22}\rho_{22}
\end{array}
\right),
\end{equation}
and therefore
\begin{equation}
\frac{V_1\rho V_1^*}{tr(V_1\rho V_1^*)}=\left(
\begin{array}{cc}
1 & 0\\
0 & 0
\end{array}
\right) ,\sp \frac{V_2\rho V_2^*}{tr(V_2\rho V_2^*)}=\left(
\begin{array}{cc}
1 & 0\\
0 & 0
\end{array}
\right)
\end{equation}
\begin{equation}
\frac{V_3\rho V_3^*}{tr(V_3\rho V_3^*)}=\left(
\begin{array}{cc}
0 & 0\\
0 & 1
\end{array}
\right)  ,\sp \frac{V_4\rho V_4^*}{tr(V_4\rho V_4^*)}=\left(
\begin{array}{cc}
0 & 0\\
0 & 1
\end{array}
\right)
\end{equation}
that is, the above values do not depend on  $\rho$.

Define
\begin{equation}
\rho_x=\left(
\begin{array}{cc}
1 & 0\\
0 & 0
\end{array}
\right) ,\sp \rho_y=\left(
\begin{array}{cc}
0 & 0\\
0 & 1
\end{array}
\right)
\end{equation}
and
\begin{equation}
\eta=\pi_1\delta_{\rho_x}+\pi_2\delta_{\rho_y}
\end{equation}
Note that the barycenter of  $\eta$ is
$$\rho_{\eta}=\pi_1\rho_x+\pi_2\rho_y=\pi_1\left(
\begin{array}{cc}
1 & 0\\
0 & 0
\end{array}
\right)+\pi_2\left(
\begin{array}{cc}
0 & 0\\
0 & 1
\end{array}
\right)=\left(
\begin{array}{cc}
\pi_1 & 0\\
0 & \pi_2
\end{array}
\right)$$
One can show directly that $\mathcal{V}(\eta)
=\eta$ (see \cite{BLLT1}). Define
\begin{equation}
\rho_1=\left(
\begin{array}{cc}
1 & 0\\
0 & 0
\end{array}
\right) ,\sp \rho_2=\left(
\begin{array}{cc}
0 & 0\\
0 & 1
\end{array}
\right)
\end{equation}
and also
\begin{equation}
\eta=\pi_1\delta_{\rho_1}+\pi_2\delta_{\rho_2}
\end{equation}
Note that the barycenter of $\eta$ is
$$\rho_{\eta}=\pi_1\rho_1+\pi_2\rho_2=\pi_1\left(
\begin{array}{cc}
1 & 0\\
0 & 0
\end{array}
\right)+\pi_2\left(
\begin{array}{cc}
0 & 0\\
0 & 1
\end{array}
\right)=\left(
\begin{array}{cc}
\pi_1 & 0\\
0 & \pi_2
\end{array}
\right)$$
From this it will also follow  that  $\mathcal{V}\eta=\eta$ \cite{BLLT1}. We will show that the entropy of such $\eta$ is $\log(2)-\pi_1\log(\pi_1)-\pi_2\log(\pi_2)$. Remember that
$$\int \log{\Big(\sum_i \frac{f\circ F_i}{f}\Big)}d\mu$$
\begin{equation}
=\int \log{\Big(\sum_i f\Big(\frac{V_i\rho V_i^*}{tr(V_i\rho V_i^*)}\Big)\Big)}d\mu-\sum_i\int tr(W_i\rho W_i^*)\log f\Big(\frac{V_i\rho V_i^*}{tr(V_i\rho V_i^*)}\Big)d\mu.
\end{equation}
For such choice of  $V_i$ take
\begin{equation}
c_i=f\Big(\frac{V_i\rho V_i^*}{tr(V_i\rho V_i^*)}\Big), \sp i=1,\dots, 4
\end{equation}
Note that
\begin{equation}\label{cond_int_set09}
c_1=c_2, \sp c_3=c_4
\end{equation}

Then we can write
\begin{equation}\label{qm_inv1}
\int \log{\Big(\sum_i \frac{f\circ F_i}{f}\Big)}d\eta
=\int \log{\Big(\sum_i c_i\Big)}d\eta-\sum_i\int tr(W_i\rho W_i^*)\log c_i d\eta.
\end{equation}

Therefore

$$\int \log{\Big(\sum_i \frac{f\circ F_i}{f}\Big)}d\eta
%$$=\pi_1\log{\Big(\sum_i f(\frac{V_i\rho_1 V_i^*}{tr(V_i\rho_1 V_i^*)}) %\Big)}+\pi_2\log{\Big(\sum_i f(\frac{V_i\rho_2 V_i^*}{tr(V_i\rho_2 V_i^*)}) \Big)}$$
=\pi_1\log{\Big(\sum_i c_i \Big)}+\pi_2\log{\Big(\sum_i c_i \Big)}$$
$$-\sum_i \Big[ tr(V_i\rho_1 V_i^*)\pi_1\log{(c_i)}+tr(V_i\rho_2 V_i^*)\pi_2\log{(c_i)}\Big]$$
$$=\pi_1\log{(2c_1+2c_3)}+\pi_2\log{(2c_1+2c_3)}$$
$$-\sum_i \Big[ tr(V_i\rho_1 V_i^*)\pi_1\log{(c_i)}+tr(V_i\rho_2 V_i^*)\pi_2\log{(c_i)}\Big]$$
$$=\log{(2(c_1+c_3))}-\sum_i \Big[ tr(V_i\rho_1 V_i^*)\pi_1\log{(c_i)}+tr(V_i\rho_2 V_i^*)\pi_2\log{(c_i)}\Big]$$
$$=\log{(2(c_1+c_3))}$$
$$-\Big[\pi_1\Big(p_{11}^2\log(c_1)+p_{12}p_{21}\log(c_2)+p_{21}p_{11}\log(c_3)+p_{22}p_{21}\log(c_4)\Big)$$
$$+\pi_2\Big(p_{11}p_{12}\log(c_1)+p_{12}p_{22}\log(c_2)+p_{21}p_{12}\log(c_3)+p_{22}^2\log(c_4)\Big) \Big]$$
$$=\log{(2(c_1+c_3))}$$
$$-\Big[p_{11}\log(c_1)(\pi_1p_{11}+\pi_2p_{12})+p_{12}\log(c_2)(\pi_1p_{21}+\pi_2p_{22})$$
$$+p_{21}\log(c_3)(\pi_1p_{11}+\pi_2p_{12})+p_{22}\log(c_4)(\pi_1p_{21}+\pi_2p_{22})\Big]$$
$$=\log{(2(c_1+c_3))}$$
$$-\Big[\pi_1p_{11}\log(c_1)+\pi_2p_{12}\log(c_2)+\pi_1p_{21}\log(c_3)+\pi_2p_{22}\log(c_4)\Big]$$
$$=\log{(2(c_1+c_3))}-(\pi_1\log(c_1)+\pi_2\log(c_3))$$
Finally,
\begin{equation}\label{eq_e_0210}
\int \log{\Big(\sum_i \frac{f\circ F_i}{f}\Big)}d\eta=\log{(2(c_1+c_3))}-(\pi_1\log(c_1)+\pi_2\log(c_3)).
\end{equation}

Now we will use Lagrange multipliers. Define $b:\mathbb{R}_+^2\to \mathbb{R}$, where $\mathbb{R}_+^2$ is the set of positive coordinates, by

$$b(x,y)=\log{(2(x+y))}-(\pi_1\log(x)+\pi_2\log(y))$$

We impose the restriction
$$x+y=a$$
 for fixed $a>0$. We will get bellow   the critical point of  $b$ under such restriction. After that we consider a general $a>0$.

Define
$$g(x,y)=x+y-a$$
and
$$\Gamma(x,y,\lambda)=b+\lambda g$$
Then, $\nabla\Gamma=0$ implies
\begin{equation}
\frac{1}{x+y}-\frac{\pi_1}{x}+\lambda=0
\end{equation}
\begin{equation}
\frac{1}{x+y}-\frac{\pi_2}{y}+\lambda=0
\end{equation}
\begin{equation}
x+y=a
\end{equation}
from which follows
\begin{equation}
x=\pi_1 a,\sp y=\pi_2 a.
\end{equation}
Therefore,
\begin{equation}
c_1=c_2=\pi_1 a, \sp c_3=c_4=\pi_2 a
\end{equation}
From (\ref{eq_e_0210}) we get
$$\int \log{\Big(\sum_i \frac{f\circ F_i}{f}\Big)}d\eta=\log{(2(\pi_1a+\pi_2a))}-(\pi_1\log(\pi_1a)+\pi_2\log(\pi_2a))$$
$$=\log(2a)-\pi_1\log(\pi_1a)-\pi_2\log(\pi_2a)$$
$$=\log(2)+\log(a)-\pi_1\log(\pi_1)-\pi_1\log(a)-\pi_2\log(\pi_2)-\pi_2\log(a)$$
\begin{equation}\label{res_021009_res}
=\log(2)-\pi_1\log(\pi_1)-\pi_2\log(\pi_2)
\end{equation}

This value of entropy is different from the value computed in the same example of QIFS
in \cite{BLLT1}, \cite{BLLT2} which is $-\sum_{i,j}\pi_i p_{ji}\log{p_{ji}}$ (Example 7 in section 11 \cite{BLLT1}).
\qee
\end{examp}

Given the expression
$$h_0(\nu):=\inf_{f \in \mathbb{B}^{+}}\int \log(\sum_{i=1}^{k} \frac{ f \circ F_{i}}{f}) d\nu,$$
for a fixed probability $\nu$,  which is invariant by the shift
acting on the space $\Omega$, a natural question is to identify the $f$
which realizes the infimum above.

We will describe below the analysis of the classical case (in the sense of Stochastic Processes, and not QSP). Our purpose is to explain why the definition presented above is a natural generalization of the setting for Markov Processes.  In the case the probability
$\nu$ comes from a Markov Process this will be now derived.

Let $\Omega=I_m^{\mathbb{N}}$, where $I_m=\{1,\dots, m\}$, and let
$\mathcal{C}=\{C_\iota:\iota\in \cup_{n\in\mathbb{N}}I_m^n\}$ the
collection of cylinder sets in $\Omega$, where
$$C_\iota:=\{\omega\in I_k^\mathbb{N}: w(j)=i_j, j=1,\dots, r, \iota=(i_1,\dots, i_r)\in I_m^r\}$$
and denote by $\sigma(\mathcal{C})$ the $\sigma$-algebra generated
by the cylinders in $\Omega$.

\bigskip Let $(P,\pi)$ be a Markov chain, so that $P=(p_{ij})$ is a
matrix of order $n$, with $p_{ij}\geq 0$, $\sum_j p_{ij}=1$ (row
stochastic), and $\pi=(\pi_1,\dots, \pi_n)$ is the left
eigenvector with eigenvalue 1. So $\pi P=\pi$, that is, $\sum_i
\pi_ip_{ij}=\pi_j $.

\bigskip

Associated to the matrix $P$ we have the following measure.

\begin{defi}
The {\bf Markov measure} (associated to the chain $(P,\pi)$) of a
cylinder is defined as
\begin{equation}
\nu(C_\iota):=\pi_{i_1}p_{i_{1}i_{2}}p_{i_{2}i_{3}}\cdots
p_{i_{r-1}i_r}
\end{equation}
\end{defi}

\qee

We  are interested in the following problem: find the infimum $f$
in
\begin{equation}\label{h0_eq}
h_0(\nu):=\inf_{f \in \mathbb{B}^{+}}\int \log(\sum_{i=1}^{k}
\frac{ f \circ F_{i}}{f}) d\nu
\end{equation}
for such $\nu$ defined above.

\qee

We use the notation $\ov{ij}$ to denote the cylinder set in
$I_m^{\mathbb{N}}$ which consists of the set of sequences
$(w_1,w_2,\dots)$ such that $w_1=i$ and $w_2=j$. Denote by
$1_{\ov{ij}}$ the indicator function of $\ov{ij}$. To simplify,
suppose $m=2$ so the alphabet considered contains only two
symbols, denoted by $1$ and $2$. Define the following function
$f:I_2^{\mathbb{N}}\to\mathbb{R}^+$,
\begin{equation}
f(x)=\sum_{i,j=1}^2 a_{ij} 1_{\ov{ij}}(x)
\end{equation}
where $a_{ij}\in\mathbb{R}^+$. That is, $f$ is a simple function,
constant on $\ov{ij}$. In this form, $\log f=\sum_{i,j} \log
a_{ij}1_{\ov{ij}}$.

\bigskip

Let us suppose that $F_i:I_m^{\mathbb{N}}\to I_m^{\mathbb{N}}$ is
the mapping $F_i(w_1,w_2,\dots)=(i,w_1,w_2,\dots)$. If $\nu$ is a
Markov measure, we have
\begin{equation}
\int_{I_2^{\mathbb{N}}}\log
fd\nu=\int_{I_2^{\mathbb{N}}}\sum_{i,j=1}^2
\log{(a_{ij})}1_{\ov{ij}} d\nu=\sum_{i,j=1}^2\log{(a_{ij})}
\nu({\ov{ij}})=\sum_{i,j=1}^2 \pi_{i}p_{ij}\log a_{ij}
\end{equation}
Also, we have for $w=(i,j,\dots)$,
\begin{equation}\label{e_aux1a}
f\circ F_l(w)=\sum_{i,j}a_{ij}1_{\ov{ij}}(F_l(w))=a_{li}
\end{equation}
To see that, note that by the expression above we have a sum of
terms such that $(i,j)=(l,i)$, therefore $a_{ij}=a_{li}$.

\bigskip

Then
$$\int \log(\sum_{i=1}^{2} \frac{ f \circ F_{i}}{f}) d\nu=\int \log (\sum_{l=1}^{2} f \circ F_{l})d\nu-\int \log f d\nu$$
\begin{equation}\label{um_v_sup}
=\int \log (\sum_{l=1}^{2} f \circ F_{l})d\nu-\sum_{i,j=1}^2
\pi_{i}p_{ij}\log a_{ij}
\end{equation}
Note that for any $w\in I_m^{\mathbb{N}}$, $w=(1,\dots)$ or
$w=(2,\dots)$. Then, by (\ref{e_aux1a}) we get
\begin{equation}\label{soma_pec}
\sum_{l=1}^2 f \circ F_l(w)=\left\{\begin{array}{ll}
a_{11}+a_{21} & \textrm{ se $w=(1,\dots)$}\\
a_{12}+a_{22} & \textrm{ se $w=(2,\dots)$}
\end{array} \right.
\end{equation}

Now fix $a_{ij}=p_{ji}$, where $p_{ij}$ are the entries of the row
stochastic matrix $P$ initially fixed. Then we get
$a_{11}+a_{21}=p_{11}+p_{12}=1$ e $a_{12}+a_{22}=p_{21}+p_{22}=1$.
Therefore for such choice of $a_{ij}$ and for any $w \in
I_m^{\mathbb{N}}$, the sum (\ref{soma_pec}) equals 1. So, by
(\ref{um_v_sup}), we get
\begin{equation}
\int \log(\sum_{i=1}^{2} \frac{ f \circ F_{i}}{f})
d\nu=-\sum_{i,j=1}^2\pi_{i}p_{ij}\log p_{ij}=H(P)
\end{equation}
Therefore,
\begin{equation}\label{um_lado001}
\inf_{f\in\mathbb{B}^+}\int \log(\sum_{i=1}^{2} \frac{ f \circ
F_{i}}{f}) d\nu\leq H(P)
\end{equation}

\qee

Now note that any positive function $f$ can be written as
$$f(w)=\sum_{i,j=1}^2 a_{ij} p_{ji}\,1_{\ov{ji}}(w)$$
Define
$$u(w):=\sum_{i,j=1}^2 a_{ij} \,1_{\ov{ji}}(w)$$
and
$$g(w):=\sum_{i,j=1}^2 p_{ji} \,1_{\ov{ji}}(w)$$
We have
$$\int_{I_2^{\mathbb{N}}}\log fd\nu=\int_{I_2^{\mathbb{N}}}\sum_{i,j=1}^2 \log{(a_{ij}p_{ji})}1_{\ov{ji}} d\nu=\sum_{i,j=1}^2\log{(a_{ij}p_{ji})} \nu({\ov{ji}})$$
\begin{equation}\label{ue_11}
=\sum_{i,j=1}^2 \pi_{j}p_{ji}\log (a_{ij}p_{ji})=\sum_{i,j=1}^2
\pi_{j}p_{ji}\log (a_{ij})+\sum_{i,j=1}^2
\pi_{j}p_{ji}\log(p_{ji})
\end{equation}
If $w=(i,j,\dots)$, then $f\circ F_l(w)=a_{li}p_{il}$ and so
$$\sum_l f\circ F_l=\sum_l a_{li}p_{il}$$
We write
\begin{equation}
\mathcal{L}_g(u)(w)=\sum_{l} f \circ
F_{l}(w)=\sum_l\sum_{i,j}a_{ij}p_{ji}1_{ij}(F_l(w))
\end{equation}
We also have the following:

\begin{lem}
\begin{equation}\label{ue_13}
\int\mathcal{L}_g(\log u) d\nu=\int\log u d\nu
\end{equation}
\end{lem}
{\bf Proof} We have
\begin{equation}\label{prec_dtb1}
\int\log u d\nu=\int\sum_{i,j}\log
(a_{ij})1_{ji}d\nu=\sum_{i,j}\log
(a_{ij})\nu(\ov{ji})=\sum_{i,j}\log (a_{ij})\pi_j p_{ji}
\end{equation}
And also
$$\int\mathcal{L}_g(\log u) d\nu=\int \sum_l\sum_{i,j}\log (a_{ij})p_{ji}1_{\ov{ij}}(F_l(w)) d\nu$$
$$=\sum_{i,j}\log (a_{ij})p_{ji}\sum_l \int 1_{\ov{ij}}(F_l(w)) d\nu$$
\begin{equation}
=\sum_{i,j}\log (a_{ij})p_{ji}\sum_l \nu(\ov{lj})=\sum_{i,j}\log
(a_{ij})p_{ji}(\pi_1 p_{1j}+\pi_2 p_{2j})=\sum_{i,j}\log
(a_{ij})\pi_j p_{ji}
\end{equation}
So,
\begin{equation}
\int\mathcal{L}_g(\log u) d\nu=\int\log u d\nu
\end{equation}
\qed Then, by using (\ref{ue_11}), (\ref{ue_13}) and
(\ref{prec_dtb1}),
$$\int \log(\sum_{i=1}^{k} \frac{ f \circ F_{i}}{f}) d\nu =  \int \log(\sum_{l=1}^{2} f \circ F_{l})d\nu-\int \log f d\nu$$
$$=\int \log(\sum_{l=1}^{2} f \circ F_{l})d\nu-\Big(\sum_{i,j=1}^2 \pi_{j}p_{ji}\log (a_{ij})+\sum_{i,j=1}^2 \pi_{i}p_{ij}\log(p_{ij})\Big)$$
\begin{equation}
=\int \log\,({\cal L}_{g}(u)) \,d \nu - \int \log u \,d \nu + H(P)
\end{equation}
\begin{equation}\label{cxv1}
=\int \log\,({\cal L}_{g}(u)) \,d \nu - \int {\cal L}_{g}(\log u)
\,d \nu + H(P)
\end{equation}
We would like to show that
\begin{equation}\label{cxv2}
\int \log\,({\cal L}_{g}(u)) \,d \nu - \int {\cal L}_{g}(\log u)
\,d \nu\geq 0
\end{equation}
This follows immediately if we show that for $w=(i,j,\dots)$,
\begin{equation}
\log\,({\cal L}_{g}(u)) \,(w) \geq {\cal L}_{g}(\log u) \,(w)
\end{equation}
The last expression follows from convexity. Indeed, to prove the
above inequality, it is enough to show that for any
$w=(i,j,\dots)$, we have
\begin{equation}
\log\Big( \sum_l a_{li}p_{il}\Big)\geq \sum_l p_{il}\log a_{li}
\end{equation}
And such inequality is true, because the $p_{il}$ are positive
numbers with $\sum_l p_{il}=1$, for any $i$, and the function
 $\log$ is concave.

\bigskip

Therefore we conclude from (\ref{cxv1}) and (\ref{cxv2}) that
\begin{equation}\label{um_lado002}
\int \log(\sum_{i=1}^{k} \frac{ f \circ F_{i}}{f}) d\nu\geq  H(P)
\end{equation}

\bigskip

{\bf Conclusion} By (\ref{um_lado001}) and (\ref{um_lado002}) we
conclude that if $\nu$ is a Markov measure associated to a
stochastic matrix $P$, then
\begin{equation}
\inf_{f\in\mathbb{B}^+}\int \log(\sum_{i=1}^{2} \frac{ f \circ
F_{i}}{f}) d\nu= H(P),
\end{equation}
and the function
$f$ such that
\begin{equation}
f(x)=\sum_{i,j=1}^2 p_{ij} 1_{\ov{ij}}(x)
\end{equation}
realizes the infimum.

\qee

We conclude this section by stating the variational problem of
pressure for our setting. We consider the the set of $V_i$, $i=1,2,\dots,k$ fixed, and we consider a variable set of $W_i$, $i=1,2,\dots,k$. In the normalized case, the different possible choices of $p_i, i=1,2,\dots,k$, (which means different choices of $W_i, i=1,2,\dots,k$) play here the role of the different Jacobians  of possible invariant probabilities (see \cite{Man} II.1, and \cite{lopes}) in Thermodynamic Formalism. In some sense the  probabilities $\mu$ can be identified with the Jacobians (this is true at least for Gibbs probabilities of H\"older potentials \cite{par}). The set of Gibbs probabilities for H\"older potentials is dense in the set of invariant probabilities \cite{lop1}.

Let $H:\mathcal{M}_N\to\mathcal{M}_N$ be
a hermitian operator. We have the following problem. Define
$F_0:\mathcal{M}_F\to\mathbb{R}$,
$$F_0(\mu):=h_0(\mu)-\frac{1}{T}tr(H\rho_{\mu})=\inf_{f \in \mathbb{B}^{+}}\int \log(\sum_{i=1}^{k} \frac{ f \circ F_{i}}{f}) d\mu -\frac{1}{T}tr(H\rho_{\mu}),$$
where $\rho_{\mu}$ is the barycenter of $\mu$, that is, the unique
$\rho\in \mathcal{M}_N$ such that
$$l(\rho)=\int_{\mathcal{M}_N} l d\mu,$$
for all $l\in V^*$. Then, in order to find the associated Gibbs
state we have   to find $\hat{\mu}\in\mathcal{M}_F$ such that
$$F_0 (\hat{\mu})=\sup_{\mu\in\mathcal{M}_F}F_0(\mu).$$

We consider above each $\mu$ which is associated to a possible set of $W_i$.
\qee

\end{document}